\input amstex
%

\def\next{AMS-SEKR}\ifx\styname\next \endinput\fi
\catcode`\@=11
\def\styname{AMS-SEKR}
\def\styversion{2.0}
{\W@{}\W@{\styname.STY - Version \styversion}\W@{}}
\hyphenation{acad-e-my acad-e-mies af-ter-thought anom-aly anom-alies
an-ti-deriv-a-tive an-tin-o-my an-tin-o-mies apoth-e-o-ses apoth-e-o-sis
ap-pen-dix ar-che-typ-al as-sign-a-ble as-sist-ant-ship as-ymp-tot-ic
asyn-chro-nous at-trib-uted at-trib-ut-able bank-rupt bank-rupt-cy
bi-dif-fer-en-tial blue-print busier busiest cat-a-stroph-ic
cat-a-stroph-i-cally con-gress cross-hatched data-base de-fin-i-tive
de-riv-a-tive dis-trib-ute dri-ver dri-vers eco-nom-ics econ-o-mist
elit-ist equi-vari-ant ex-quis-ite ex-tra-or-di-nary flow-chart
for-mi-da-ble forth-right friv-o-lous ge-o-des-ic ge-o-det-ic geo-met-ric
griev-ance griev-ous griev-ous-ly hexa-dec-i-mal ho-lo-no-my ho-mo-thetic
ideals idio-syn-crasy in-fin-ite-ly in-fin-i-tes-i-mal ir-rev-o-ca-ble
key-stroke lam-en-ta-ble light-weight mal-a-prop-ism man-u-script
mar-gin-al meta-bol-ic me-tab-o-lism meta-lan-guage me-trop-o-lis
met-ro-pol-i-tan mi-nut-est mol-e-cule mono-chrome mono-pole mo-nop-oly
mono-spline mo-not-o-nous mul-ti-fac-eted mul-ti-plic-able non-euclid-ean
non-iso-mor-phic non-smooth par-a-digm par-a-bol-ic pa-rab-o-loid
pa-ram-e-trize para-mount pen-ta-gon phe-nom-e-non post-script pre-am-ble
pro-ce-dur-al pro-hib-i-tive pro-hib-i-tive-ly pseu-do-dif-fer-en-tial
pseu-do-fi-nite pseu-do-nym qua-drat-ics quad-ra-ture qua-si-smooth
qua-si-sta-tion-ary qua-si-tri-an-gu-lar quin-tes-sence quin-tes-sen-tial
re-arrange-ment rec-tan-gle ret-ri-bu-tion retro-fit retro-fit-ted
right-eous right-eous-ness ro-bot ro-bot-ics sched-ul-ing se-mes-ter
semi-def-i-nite semi-ho-mo-thet-ic set-up se-vere-ly side-step sov-er-eign
spe-cious spher-oid spher-oid-al star-tling star-tling-ly
sta-tis-tics sto-chas-tic straight-est strange-ness strat-a-gem strong-hold
sum-ma-ble symp-to-matic syn-chro-nous topo-graph-i-cal tra-vers-a-ble
tra-ver-sal tra-ver-sals treach-ery turn-around un-at-tached un-err-ing-ly
white-space wide-spread wing-spread wretch-ed wretch-ed-ly Brown-ian
Eng-lish Euler-ian Feb-ru-ary Gauss-ian Grothen-dieck Hamil-ton-ian
Her-mit-ian Jan-u-ary Japan-ese Kor-te-weg Le-gendre Lip-schitz
Lip-schitz-ian Mar-kov-ian Noe-ther-ian No-vem-ber Rie-mann-ian
Schwarz-schild Sep-tem-ber
form per-iods Uni-ver-si-ty cri-ti-sism for-ma-lism}
\Invalid@\nofrills
\Invalid@\usualspace
\newif\ifnofrills@
\def\nofrills@#1#2{\relaxnext@
  \DN@{\ifx\next\nofrills
    \nofrills@true\let#2\relax\DN@\nofrills{\nextii@}%
  \else
    \nofrills@false\def#2{#1}\let\next@\nextii@\fi
\next@}}
\def\usualspace@#1{\ifnofrills@\def\usualspace{#1}\fi}
\def\addto#1#2{\csname \expandafter\eat@\string#1@\endcsname
  \expandafter{\the\csname \expandafter\eat@\string#1@\endcsname#2}}
\newdimen\bigsize@
\def\big@#1#2{{\hbox{$\left#2\vcenter to#1\bigsize@{}%
  \right.\nulldelimiterspace\z@\m@th$}}}
\def\big{\big@\@ne}
\def\Big{\big@{1.5}}
\def\bigg{\big@\tw@}
\def\Bigg{\big@{2.5}}
\def\raggedcenter@{\leftskip\z@ plus.4\hsize \rightskip\leftskip
 \parfillskip\z@ \parindent\z@ \spaceskip.3333em \xspaceskip.5em
 \pretolerance9999\tolerance9999 \exhyphenpenalty\@M
 \hyphenpenalty\@M \let\\\linebreak}
\def\upperspecialchars{\def\ss{SS}\let\i=I\let\j=J\let\ae\AE\let\oe\OE
  \let\o\O\let\aa\AA\let\l\L}
\def\uppercasetext@#1{%
  {\spaceskip1.2\fontdimen2\the\font plus1.2\fontdimen3\the\font
   \upperspecialchars\uctext@#1$\m@th\aftergroup\eat@$}}
\def\uctext@#1$#2${\endash@#1-\endash@$#2$\uctext@}
\def\endash@#1-#2\endash@{\uppercase{#1}\if\notempty{#2}--\endash@#2\endash@\fi}
\def\runaway@#1{\DN@{#1}\ifx\envir@\next@
  \Err@{You seem to have a missing or misspelled \string\end#1 ...}%
  \let\envir@\empty\fi}
\newif\iftemp@
\def\notempty#1{TT\fi\def\test@{#1}\ifx\test@\empty\temp@false
  \else\temp@true\fi \iftemp@}
\font@\tensmc=cmcsc10
\font@\sevenex=cmex7
\font@\sevenit=cmti7
\font@\eightrm=cmr8 
\font@\sixrm=cmr6 
\font@\eighti=cmmi8     \skewchar\eighti='177 
\font@\sixi=cmmi6       \skewchar\sixi='177   
\font@\eightsy=cmsy8    \skewchar\eightsy='60 
\font@\sixsy=cmsy6      \skewchar\sixsy='60   
\font@\eightex=cmex8
\font@\eightbf=cmbx8 
\font@\sixbf=cmbx6   
\font@\eightit=cmti8 
\font@\eightsl=cmsl8 
\font@\eightsmc=cmcsc8
\font@\eighttt=cmtt8 


\loadmsam
\loadmsbm
\loadeufm
\UseAMSsymbols
\newtoks\tenpoint@
\def\tenpoint{\normalbaselineskip12\p@
 \abovedisplayskip12\p@ plus3\p@ minus9\p@
 \belowdisplayskip\abovedisplayskip
 \abovedisplayshortskip\z@ plus3\p@
 \belowdisplayshortskip7\p@ plus3\p@ minus4\p@
 \textonlyfont@\rm\tenrm \textonlyfont@\it\tenit
 \textonlyfont@\sl\tensl \textonlyfont@\bf\tenbf
 \textonlyfont@\smc\tensmc \textonlyfont@\tt\tentt
 \textonlyfont@\bsmc\tenbsmc
 \ifsyntax@ \def\big##1{{\hbox{$\left##1\right.$}}}%
  \let\Big\big \let\bigg\big \let\Bigg\big
 \else
  \textfont\z@=\tenrm  \scriptfont\z@=\sevenrm  \scriptscriptfont\z@=\fiverm
  \textfont\@ne=\teni  \scriptfont\@ne=\seveni  \scriptscriptfont\@ne=\fivei
  \textfont\tw@=\tensy \scriptfont\tw@=\sevensy \scriptscriptfont\tw@=\fivesy
  \textfont\thr@@=\tenex \scriptfont\thr@@=\sevenex
        \scriptscriptfont\thr@@=\sevenex
  \textfont\itfam=\tenit \scriptfont\itfam=\sevenit
        \scriptscriptfont\itfam=\sevenit
  \textfont\bffam=\tenbf \scriptfont\bffam=\sevenbf
        \scriptscriptfont\bffam=\fivebf
  \setbox\strutbox\hbox{\vrule height8.5\p@ depth3.5\p@ width\z@}%
  \setbox\strutbox@\hbox{\lower.5\normallineskiplimit\vbox{%
        \kern-\normallineskiplimit\copy\strutbox}}%
 \setbox\z@\vbox{\hbox{$($}\kern\z@}\bigsize@=1.2\ht\z@
 \fi
 \normalbaselines\rm\ex@.2326ex\jot3\ex@\the\tenpoint@}
\newtoks\eightpoint@
\def\eightpoint{\normalbaselineskip10\p@
 \abovedisplayskip10\p@ plus2.4\p@ minus7.2\p@
 \belowdisplayskip\abovedisplayskip
 \abovedisplayshortskip\z@ plus2.4\p@
 \belowdisplayshortskip5.6\p@ plus2.4\p@ minus3.2\p@
 \textonlyfont@\rm\eightrm \textonlyfont@\it\eightit
 \textonlyfont@\sl\eightsl \textonlyfont@\bf\eightbf
 \textonlyfont@\smc\eightsmc \textonlyfont@\tt\eighttt
 \textonlyfont@\bsmc\eightbsmc
 \ifsyntax@\def\big##1{{\hbox{$\left##1\right.$}}}%
  \let\Big\big \let\bigg\big \let\Bigg\big
 \else
  \textfont\z@=\eightrm \scriptfont\z@=\sixrm \scriptscriptfont\z@=\fiverm
  \textfont\@ne=\eighti \scriptfont\@ne=\sixi \scriptscriptfont\@ne=\fivei
  \textfont\tw@=\eightsy \scriptfont\tw@=\sixsy \scriptscriptfont\tw@=\fivesy
  \textfont\thr@@=\eightex \scriptfont\thr@@=\sevenex
   \scriptscriptfont\thr@@=\sevenex
  \textfont\itfam=\eightit \scriptfont\itfam=\sevenit
   \scriptscriptfont\itfam=\sevenit
  \textfont\bffam=\eightbf \scriptfont\bffam=\sixbf
   \scriptscriptfont\bffam=\fivebf
 \setbox\strutbox\hbox{\vrule height7\p@ depth3\p@ width\z@}%
 \setbox\strutbox@\hbox{\raise.5\normallineskiplimit\vbox{%
   \kern-\normallineskiplimit\copy\strutbox}}%
 \setbox\z@\vbox{\hbox{$($}\kern\z@}\bigsize@=1.2\ht\z@
 \fi
 \normalbaselines\eightrm\ex@.2326ex\jot3\ex@\the\eightpoint@}

\font@\twelverm=cmr10 scaled\magstep1
\font@\twelveit=cmti10 scaled\magstep1
\font@\twelvesl=cmsl10 scaled\magstep1
\font@\twelvesmc=cmcsc10 scaled\magstep1
\font@\twelvett=cmtt10 scaled\magstep1
\font@\twelvebf=cmbx10 scaled\magstep1
\font@\twelvei=cmmi10 scaled\magstep1
\font@\twelvesy=cmsy10 scaled\magstep1
\font@\twelveex=cmex10 scaled\magstep1
\font@\twelvemsa=msam10 scaled\magstep1
\font@\twelveeufm=eufm10 scaled\magstep1
\font@\twelvemsb=msbm10 scaled\magstep1
\newtoks\twelvepoint@
\def\twelvepoint{\normalbaselineskip15\p@
 \abovedisplayskip15\p@ plus3.6\p@ minus10.8\p@
 \belowdisplayskip\abovedisplayskip
 \abovedisplayshortskip\z@ plus3.6\p@
 \belowdisplayshortskip8.4\p@ plus3.6\p@ minus4.8\p@
 \textonlyfont@\rm\twelverm \textonlyfont@\it\twelveit
 \textonlyfont@\sl\twelvesl \textonlyfont@\bf\twelvebf
 \textonlyfont@\smc\twelvesmc \textonlyfont@\tt\twelvett
 \textonlyfont@\bsmc\twelvebsmc
 \ifsyntax@ \def\big##1{{\hbox{$\left##1\right.$}}}%
  \let\Big\big \let\bigg\big \let\Bigg\big
 \else
  \textfont\z@=\twelverm  \scriptfont\z@=\tenrm  \scriptscriptfont\z@=\sevenrm
  \textfont\@ne=\twelvei  \scriptfont\@ne=\teni  \scriptscriptfont\@ne=\seveni
  \textfont\tw@=\twelvesy \scriptfont\tw@=\tensy \scriptscriptfont\tw@=\sevensy
  \textfont\thr@@=\twelveex \scriptfont\thr@@=\tenex
        \scriptscriptfont\thr@@=\tenex
  \textfont\itfam=\twelveit \scriptfont\itfam=\tenit
        \scriptscriptfont\itfam=\tenit
  \textfont\bffam=\twelvebf \scriptfont\bffam=\tenbf
        \scriptscriptfont\bffam=\sevenbf
  \setbox\strutbox\hbox{\vrule height10.2\p@ depth4.2\p@ width\z@}%
  \setbox\strutbox@\hbox{\lower.6\normallineskiplimit\vbox{%
        \kern-\normallineskiplimit\copy\strutbox}}%
 \setbox\z@\vbox{\hbox{$($}\kern\z@}\bigsize@=1.4\ht\z@
 \fi
 \normalbaselines\rm\ex@.2326ex\jot3.6\ex@\the\twelvepoint@}

\def\headfonts{\twelvepoint\bf}

\font@\fourteenrm=cmr10 scaled\magstep2
\font@\fourteenit=cmti10 scaled\magstep2
\font@\fourteensl=cmsl10 scaled\magstep2
\font@\fourteensmc=cmcsc10 scaled\magstep2
\font@\fourteentt=cmtt10 scaled\magstep2
\font@\fourteenbf=cmbx10 scaled\magstep2
\font@\fourteeni=cmmi10 scaled\magstep2
\font@\fourteensy=cmsy10 scaled\magstep2
\font@\fourteenex=cmex10 scaled\magstep2
\font@\fourteenmsa=msam10 scaled\magstep2
\font@\fourteeneufm=eufm10 scaled\magstep2
\font@\fourteenmsb=msbm10 scaled\magstep2
\newtoks\fourteenpoint@
\def\fourteenpoint{\normalbaselineskip15\p@
 \abovedisplayskip18\p@ plus4.3\p@ minus12.9\p@
 \belowdisplayskip\abovedisplayskip
 \abovedisplayshortskip\z@ plus4.3\p@
 \belowdisplayshortskip10.1\p@ plus4.3\p@ minus5.8\p@
 \textonlyfont@\rm\fourteenrm \textonlyfont@\it\fourteenit
 \textonlyfont@\sl\fourteensl \textonlyfont@\bf\fourteenbf
 \textonlyfont@\smc\fourteensmc \textonlyfont@\tt\fourteentt
 \textonlyfont@\bsmc\fourteenbsmc
 \ifsyntax@ \def\big##1{{\hbox{$\left##1\right.$}}}%
  \let\Big\big \let\bigg\big \let\Bigg\big
 \else
  \textfont\z@=\fourteenrm  \scriptfont\z@=\twelverm  \scriptscriptfont\z@=\tenrm
  \textfont\@ne=\fourteeni  \scriptfont\@ne=\twelvei  \scriptscriptfont\@ne=\teni
  \textfont\tw@=\fourteensy \scriptfont\tw@=\twelvesy \scriptscriptfont\tw@=\tensy
  \textfont\thr@@=\fourteenex \scriptfont\thr@@=\twelveex
        \scriptscriptfont\thr@@=\twelveex
  \textfont\itfam=\fourteenit \scriptfont\itfam=\twelveit
        \scriptscriptfont\itfam=\twelveit
  \textfont\bffam=\fourteenbf \scriptfont\bffam=\twelvebf
        \scriptscriptfont\bffam=\tenbf
  \setbox\strutbox\hbox{\vrule height12.2\p@ depth5\p@ width\z@}%
  \setbox\strutbox@\hbox{\lower.72\normallineskiplimit\vbox{%
        \kern-\normallineskiplimit\copy\strutbox}}%
 \setbox\z@\vbox{\hbox{$($}\kern\z@}\bigsize@=1.7\ht\z@
 \fi
 \normalbaselines\rm\ex@.2326ex\jot4.3\ex@\the\fourteenpoint@}

\def\chapheadfonts{\fourteenpoint\bf}

\font@\seventeenrm=cmr10 scaled\magstep3
\font@\seventeenit=cmti10 scaled\magstep3
\font@\seventeensl=cmsl10 scaled\magstep3
\font@\seventeensmc=cmcsc10 scaled\magstep3
\font@\seventeentt=cmtt10 scaled\magstep3
\font@\seventeenbf=cmbx10 scaled\magstep3
\font@\seventeeni=cmmi10 scaled\magstep3
\font@\seventeensy=cmsy10 scaled\magstep3
\font@\seventeenex=cmex10 scaled\magstep3
\font@\seventeenmsa=msam10 scaled\magstep3
\font@\seventeeneufm=eufm10 scaled\magstep3
\font@\seventeenmsb=msbm10 scaled\magstep3
\newtoks\seventeenpoint@
\def\seventeenpoint{\normalbaselineskip18\p@
 \abovedisplayskip21.6\p@ plus5.2\p@ minus15.4\p@
 \belowdisplayskip\abovedisplayskip
 \abovedisplayshortskip\z@ plus5.2\p@
 \belowdisplayshortskip12.1\p@ plus5.2\p@ minus7\p@
 \textonlyfont@\rm\seventeenrm \textonlyfont@\it\seventeenit
 \textonlyfont@\sl\seventeensl \textonlyfont@\bf\seventeenbf
 \textonlyfont@\smc\seventeensmc \textonlyfont@\tt\seventeentt
 \textonlyfont@\bsmc\seventeenbsmc
 \ifsyntax@ \def\big##1{{\hbox{$\left##1\right.$}}}%
  \let\Big\big \let\bigg\big \let\Bigg\big
 \else
  \textfont\z@=\seventeenrm  \scriptfont\z@=\fourteenrm  \scriptscriptfont\z@=\twelverm
  \textfont\@ne=\seventeeni  \scriptfont\@ne=\fourteeni  \scriptscriptfont\@ne=\twelvei
  \textfont\tw@=\seventeensy \scriptfont\tw@=\fourteensy \scriptscriptfont\tw@=\twelvesy
  \textfont\thr@@=\seventeenex \scriptfont\thr@@=\fourteenex
        \scriptscriptfont\thr@@=\fourteenex
  \textfont\itfam=\seventeenit \scriptfont\itfam=\fourteenit
        \scriptscriptfont\itfam=\fourteenit
  \textfont\bffam=\seventeenbf \scriptfont\bffam=\fourteenbf
        \scriptscriptfont\bffam=\twelvebf
  \setbox\strutbox\hbox{\vrule height14.6\p@ depth6\p@ width\z@}%
  \setbox\strutbox@\hbox{\lower.86\normallineskiplimit\vbox{%
        \kern-\normallineskiplimit\copy\strutbox}}%
 \setbox\z@\vbox{\hbox{$($}\kern\z@}\bigsize@=2\ht\z@
 \fi
 \normalbaselines\rm\ex@.2326ex\jot5.2\ex@\the\seventeenpoint@}

\font@\rrrrrm=cmr10 scaled\magstep4
\font@\bigtitlefont=cmbx10 scaled\magstep4

\parindent1pc
\normallineskiplimit\p@
\newdimen\indenti \indenti=2pc
\def\pageheight#1{\vsize#1}
\def\pagewidth#1{\hsize#1%
   \captionwidth@\hsize \advance\captionwidth@-2\indenti}
\pagewidth{30pc} \pageheight{47pc}
\def\topmatter{%
 \ifx\undefined\msafam
 \else\font@\eightmsa=msam8 \font@\sixmsa=msam6
   \ifsyntax@\else \addto\tenpoint{\textfont\msafam=\tenmsa
              \scriptfont\msafam=\sevenmsa \scriptscriptfont\msafam=\fivemsa}%
     \addto\eightpoint{\textfont\msafam=\eightmsa \scriptfont\msafam=\sixmsa
              \scriptscriptfont\msafam=\fivemsa}%
   \fi
 \fi
 \ifx\undefined\msbfam
 \else\font@\eightmsb=msbm8 \font@\sixmsb=msbm6
   \ifsyntax@\else \addto\tenpoint{\textfont\msbfam=\tenmsb
         \scriptfont\msbfam=\sevenmsb \scriptscriptfont\msbfam=\fivemsb}%
     \addto\eightpoint{\textfont\msbfam=\eightmsb \scriptfont\msbfam=\sixmsb
         \scriptscriptfont\msbfam=\fivemsb}%
   \fi
 \fi
 \ifx\undefined\eufmfam
 \else \font@\eighteufm=eufm8 \font@\sixeufm=eufm6
   \ifsyntax@\else \addto\tenpoint{\textfont\eufmfam=\teneufm
       \scriptfont\eufmfam=\seveneufm \scriptscriptfont\eufmfam=\fiveeufm}%
     \addto\eightpoint{\textfont\eufmfam=\eighteufm
       \scriptfont\eufmfam=\sixeufm \scriptscriptfont\eufmfam=\fiveeufm}%
   \fi
 \fi
 \ifx\undefined\eufbfam
 \else \font@\eighteufb=eufb8 \font@\sixeufb=eufb6
   \ifsyntax@\else \addto\tenpoint{\textfont\eufbfam=\teneufb
      \scriptfont\eufbfam=\seveneufb \scriptscriptfont\eufbfam=\fiveeufb}%
    \addto\eightpoint{\textfont\eufbfam=\eighteufb
      \scriptfont\eufbfam=\sixeufb \scriptscriptfont\eufbfam=\fiveeufb}%
   \fi
 \fi
 \ifx\undefined\eusmfam
 \else \font@\eighteusm=eusm8 \font@\sixeusm=eusm6
   \ifsyntax@\else \addto\tenpoint{\textfont\eusmfam=\teneusm
       \scriptfont\eusmfam=\seveneusm \scriptscriptfont\eusmfam=\fiveeusm}%
     \addto\eightpoint{\textfont\eusmfam=\eighteusm
       \scriptfont\eusmfam=\sixeusm \scriptscriptfont\eusmfam=\fiveeusm}%
   \fi
 \fi
 \ifx\undefined\eusbfam
 \else \font@\eighteusb=eusb8 \font@\sixeusb=eusb6
   \ifsyntax@\else \addto\tenpoint{\textfont\eusbfam=\teneusb
       \scriptfont\eusbfam=\seveneusb \scriptscriptfont\eusbfam=\fiveeusb}%
     \addto\eightpoint{\textfont\eusbfam=\eighteusb
       \scriptfont\eusbfam=\sixeusb \scriptscriptfont\eusbfam=\fiveeusb}%
   \fi
 \fi
 \ifx\undefined\eurmfam
 \else \font@\eighteurm=eurm8 \font@\sixeurm=eurm6
   \ifsyntax@\else \addto\tenpoint{\textfont\eurmfam=\teneurm
       \scriptfont\eurmfam=\seveneurm \scriptscriptfont\eurmfam=\fiveeurm}%
     \addto\eightpoint{\textfont\eurmfam=\eighteurm
       \scriptfont\eurmfam=\sixeurm \scriptscriptfont\eurmfam=\fiveeurm}%
   \fi
 \fi
 \ifx\undefined\eurbfam
 \else \font@\eighteurb=eurb8 \font@\sixeurb=eurb6
   \ifsyntax@\else \addto\tenpoint{\textfont\eurbfam=\teneurb
       \scriptfont\eurbfam=\seveneurb \scriptscriptfont\eurbfam=\fiveeurb}%
    \addto\eightpoint{\textfont\eurbfam=\eighteurb
       \scriptfont\eurbfam=\sixeurb \scriptscriptfont\eurbfam=\fiveeurb}%
   \fi
 \fi
 \ifx\undefined\cmmibfam
 \else \font@\eightcmmib=cmmib8 \font@\sixcmmib=cmmib6
   \ifsyntax@\else \addto\tenpoint{\textfont\cmmibfam=\tencmmib
       \scriptfont\cmmibfam=\sevencmmib \scriptscriptfont\cmmibfam=\fivecmmib}%
    \addto\eightpoint{\textfont\cmmibfam=\eightcmmib
       \scriptfont\cmmibfam=\sixcmmib \scriptscriptfont\cmmibfam=\fivecmmib}%
   \fi
 \fi
 \ifx\undefined\cmbsyfam
 \else \font@\eightcmbsy=cmbsy8 \font@\sixcmbsy=cmbsy6
   \ifsyntax@\else \addto\tenpoint{\textfont\cmbsyfam=\tencmbsy
      \scriptfont\cmbsyfam=\sevencmbsy \scriptscriptfont\cmbsyfam=\fivecmbsy}%
    \addto\eightpoint{\textfont\cmbsyfam=\eightcmbsy
      \scriptfont\cmbsyfam=\sixcmbsy \scriptscriptfont\cmbsyfam=\fivecmbsy}%
   \fi
 \fi
 \let\topmatter\relax}
\def\chapterno@{\uppercase\expandafter{\romannumeral\chaptercount@}}
\newcount\chaptercount@
\def\chapter{\nofrills@{\afterassignment\chapterno@
                        CHAPTER \global\chaptercount@=}\chapter@
 \DNii@##1{\leavevmode\hskip-\leftskip
   \rlap{\vbox to\z@{\vss\centerline{\eightpoint
   \chapter@##1\unskip}\baselineskip2pc\null}}\hskip\leftskip
   \nofrills@false}%
 \FN@\next@}
\newbox\titlebox@

\def\title{\nofrills@{\relax}\title@%
 \DNii@##1\endtitle{\global\setbox\titlebox@\vtop{\tenpoint\bf
 \raggedcenter@\ignorespaces
 \baselineskip1.3\baselineskip\title@{##1}\endgraf}%
 \ifmonograph@ \edef\next{\the\leftheadtoks}\ifx\next\empty
    \leftheadtext{##1}\fi
 \fi
 \edef\next{\the\rightheadtoks}\ifx\next\empty \rightheadtext{##1}\fi
 }\FN@\next@}
\newbox\authorbox@
\def\author#1\endauthor{\global\setbox\authorbox@
 \vbox{\tenpoint\smc\raggedcenter@\ignorespaces
 #1\endgraf}\relaxnext@ \edef\next{\the\leftheadtoks}%
 \ifx\next\empty\leftheadtext{#1}\fi}
\newbox\affilbox@
\def\affil#1\endaffil{\global\setbox\affilbox@
 \vbox{\tenpoint\raggedcenter@\ignorespaces#1\endgraf}}
\newcount\addresscount@
\addresscount@\z@
\def\address#1\endaddress{\global\advance\addresscount@\@ne
  \expandafter\gdef\csname address\number\addresscount@\endcsname
  {\vskip12\p@ minus6\p@\noindent\eightpoint\smc\ignorespaces#1\par}}
\def\email{\nofrills@{\eightpoint{\it E-mail\/}:\enspace}\email@
  \DNii@##1\endemail{%
  \expandafter\gdef\csname email\number\addresscount@\endcsname
  {\def\usualspace{{\it\enspace}}\smallskip\noindent\eightpoint\email@
  \ignorespaces##1\par}}%
 \FN@\next@}
\def\thedate@{}
\def\date#1\enddate{\gdef\thedate@{\tenpoint\ignorespaces#1\unskip}}
\def\thethanks@{}
\def\thanks#1\endthanks{\gdef\thethanks@{\eightpoint\ignorespaces#1.\unskip}}
\def\thekeywords@{}
\def\keywords{\nofrills@{{\it Key words and phrases.\enspace}}\keywords@
 \DNii@##1\endkeywords{\def\thekeywords@{\def\usualspace{{\it\enspace}}%
 \eightpoint\keywords@\ignorespaces##1\unskip.}}%
 \FN@\next@}
\def\thesubjclass@{}
\def\subjclass{\nofrills@{{\rm2000 {\it Mathematics Subject
   Classification\/}.\enspace}}\subjclass@
 \DNii@##1\endsubjclass{\def\thesubjclass@{\def\usualspace
  {{\rm\enspace}}\eightpoint\subjclass@\ignorespaces##1\unskip.}}%
 \FN@\next@}
\newbox\abstractbox@
\def\abstract{\nofrills@{{\smc Abstract.\enspace}}\abstract@
 \DNii@{\setbox\abstractbox@\vbox\bgroup\noindent$$\vbox\bgroup
  \def\envir@{abstract}\advance\hsize-2\indenti
  \usualspace@{{\enspace}}\eightpoint \noindent\abstract@\ignorespaces}%
 \FN@\next@}
\def\endabstract{\par\unskip\egroup$$\egroup}
\def\widestnumber#1#2{\begingroup\let\head\null\let\subhead\empty
   \let\subsubhead\subhead
   \ifx#1\head\global\setbox\tocheadbox@\hbox{#2.\enspace}%
   \else\ifx#1\subhead\global\setbox\tocsubheadbox@\hbox{#2.\enspace}%
   \else\ifx#1\key\bgroup\let\endrefitem@\egroup
        \key#2\endrefitem@\global\refindentwd\wd\keybox@
   \else\ifx#1\no\bgroup\let\endrefitem@\egroup
        \no#2\endrefitem@\global\refindentwd\wd\nobox@
   \else\ifx#1\page\global\setbox\pagesbox@\hbox{\quad\bf#2}%
   \else\ifx#1\item\setboxz@h{#2}\global\rosteritemwd\wdz@
        \global\advance\rosteritemwd by.5\parindent
   \else\message{\string\widestnumber is not defined for this option
   (\string#1)}%
\fi\fi\fi\fi\fi\fi\endgroup}
\newif\ifmonograph@
\def\Monograph{\monograph@true \let\headmark\rightheadtext
  \let\varindent@\indent \def\headfont@{\bf}\def\proclaimheadfont@{\smc}%
  \def\demofont@{\smc}}
\let\varindent@\indent

\newbox\tocheadbox@    \newbox\tocsubheadbox@
\newbox\tocbox@
\def\toc{\toc@{Contents}}
\def\newtocdefs{%
   \def \title##1\endtitle
       {\penaltyandskip@\z@\smallskipamount
        \hangindent\wd\tocheadbox@\noindent{\bf##1}}%
   \def \chapter##1{%
        Chapter \uppercase\expandafter{\romannumeral##1.\unskip}\enspace}%
   \def \specialhead##1\endspecialhead
       {\par\hangindent\wd\tocheadbox@ \noindent##1\par}%
   \def \head##1 ##2\endhead
       {\par\hangindent\wd\tocheadbox@ \noindent
        \if\notempty{##1}\hbox to\wd\tocheadbox@{\hfil##1\enspace}\fi
        ##2\par}%
   \def \subhead##1 ##2\endsubhead
       {\par\vskip-\parskip {\normalbaselines
        \advance\leftskip\wd\tocheadbox@
        \hangindent\wd\tocsubheadbox@ \noindent
        \if\notempty{##1}\hbox to\wd\tocsubheadbox@{##1\unskip\hfil}\fi
         ##2\par}}%
   \def \subsubhead##1 ##2\endsubsubhead
       {\par\vskip-\parskip {\normalbaselines
        \advance\leftskip\wd\tocheadbox@
        \hangindent\wd\tocsubheadbox@ \noindent
        \if\notempty{##1}\hbox to\wd\tocsubheadbox@{##1\unskip\hfil}\fi
        ##2\par}}}
\def\toc@#1{\relaxnext@
   \def\page##1%
       {\unskip\penalty0\null\hfil
        \rlap{\hbox to\wd\pagesbox@{\quad\hfil##1}}\hfilneg\penalty\@M}%
 \DN@{\ifx\next\nofrills\DN@\nofrills{\nextii@}%
      \else\DN@{\nextii@{{#1}}}\fi
      \next@}%
 \DNii@##1{%
\ifmonograph@\bgroup\else\setbox\tocbox@\vbox\bgroup
   \centerline{\headfont@\ignorespaces##1\unskip}\nobreak
   \vskip\belowheadskip \fi
   \setbox\tocheadbox@\hbox{0.\enspace}%
   \setbox\tocsubheadbox@\hbox{0.0.\enspace}%
   \leftskip\indenti \rightskip\leftskip
   \setbox\pagesbox@\hbox{\bf\quad000}\advance\rightskip\wd\pagesbox@
   \newtocdefs
 }%
 \FN@\next@}
\def\endtoc{\par\egroup}
\let\pretitle\relax
\let\preauthor\relax
\let\preaffil\relax
\let\predate\relax
\let\preabstract\relax
\let\prepaper\relax
\def\dedicatory #1\enddedicatory{\def\preabstract{{\medskip
  \eightpoint\it \raggedcenter@#1\endgraf}}}
\def\thetranslator@{}
\def\translator#1\endtranslator{\def\thetranslator@{\nobreak\medskip
 \line{\eightpoint\hfil Translated by \uppercase{#1}\qquad\qquad}\nobreak}}
\outer\def\endtopmatter{\runaway@{abstract}%
 \edef\next{\the\leftheadtoks}\ifx\next\empty
  \expandafter\leftheadtext\expandafter{\the\rightheadtoks}\fi
 \ifmonograph@\else
   \ifx\thesubjclass@\empty\else \makefootnote@{}{\thesubjclass@}\fi
   \ifx\thekeywords@\empty\else \makefootnote@{}{\thekeywords@}\fi
   \ifx\thethanks@\empty\else \makefootnote@{}{\thethanks@}\fi
 \fi
  \pretitle
  \ifmonograph@ \topskip7pc \else \topskip4pc \fi
  \box\titlebox@
  \topskip10pt
  \preauthor
  \ifvoid\authorbox@\else \vskip2.5pc plus1pc \unvbox\authorbox@\fi
  \preaffil
  \ifvoid\affilbox@\else \vskip1pc plus.5pc \unvbox\affilbox@\fi
  \predate
  \ifx\thedate@\empty\else \vskip1pc plus.5pc \line{\hfil\thedate@\hfil}\fi
  \preabstract
  \ifvoid\abstractbox@\else \vskip1.5pc plus.5pc \unvbox\abstractbox@ \fi
  \ifvoid\tocbox@\else\vskip1.5pc plus.5pc \unvbox\tocbox@\fi
  \prepaper
  \vskip2pc plus1pc
}
\def\document{\let\fontlist@\relax\let\alloclist@\relax
  \tenpoint}

\newskip\aboveheadskip       \aboveheadskip1.8\bigskipamount
\newdimen\belowheadskip      \belowheadskip1.8\medskipamount

\def\headfont@{\smc}
\def\penaltyandskip@#1#2{\relax\ifdim\lastskip<#2\relax\removelastskip
      \ifnum#1=\z@\else\penalty@#1\relax\fi\vskip#2%
  \else\ifnum#1=\z@\else\penalty@#1\relax\fi\fi}
\def\nobreak{\penalty\@M
  \ifvmode\def\penalty@{\let\penalty@\penalty\count@@@}%
  \everypar{\let\penalty@\penalty\everypar{}}\fi}
\let\penalty@\penalty
\def\heading#1\endheading{\head#1\endhead}

\def\specialheadfont@{\bf}
\outer\def\specialhead{\par\penaltyandskip@{-200}\aboveheadskip
  \begingroup\interlinepenalty\@M\rightskip\z@ plus\hsize \let\\\linebreak
  \specialheadfont@\noindent\ignorespaces}
\def\endspecialhead{\par\endgroup\nobreak\vskip\belowheadskip}
\let\headmark\eat@
\newskip\subheadskip       \subheadskip\medskipamount
\def\subheadfont@{\bf}
\outer\def\subhead{\nofrills@{.\enspace}\subhead@
 \DNii@##1\endsubhead{\par\penaltyandskip@{-100}\subheadskip
  \varindent@{\usualspace@{{\subheadfont@\enspace}}%
 \subheadfont@\ignorespaces##1\unskip\subhead@}\ignorespaces}%
 \FN@\next@}
\outer\def\subsubhead{\nofrills@{.\enspace}\subsubhead@
 \DNii@##1\endsubsubhead{\par\penaltyandskip@{-50}\medskipamount
      {\usualspace@{{\it\enspace}}%
  \it\ignorespaces##1\unskip\subsubhead@}\ignorespaces}%
 \FN@\next@}
\def\proclaimheadfont@{\bf}
\outer\def\proclaim{\runaway@{proclaim}\def\envir@{proclaim}%
  \nofrills@{.\enspace}\proclaim@
 \DNii@##1{\penaltyandskip@{-100}\medskipamount\varindent@
   \usualspace@{{\proclaimheadfont@\enspace}}\proclaimheadfont@
   \ignorespaces##1\unskip\proclaim@
  \sl\ignorespaces}%
 \FN@\next@}
\outer\def\endproclaim{\let\envir@\relax\par\rm
  \penaltyandskip@{55}\medskipamount}
\def\demoheadfont@{\it}
\def\demo{\runaway@{proclaim}\nofrills@{.\enspace}\demo@
     \DNii@##1{\par\penaltyandskip@\z@\medskipamount
  {\usualspace@{{\demoheadfont@\enspace}}%
  \varindent@\demoheadfont@\ignorespaces##1\unskip\demo@}\rm
  \ignorespaces}\FN@\next@}
\def\enddemo{\par\medskip}
\def\qed{\ifhmode\unskip\nobreak\fi\quad\ifmmode\square\else$\m@th\square$\fi}
\let\remark\demo
\let\endremark\enddemo
\def\definition{\runaway@{proclaim}%
  \nofrills@{.\demoheadfont@\enspace}\definition@
        \DNii@##1{\penaltyandskip@{-100}\medskipamount
        {\usualspace@{{\demoheadfont@\enspace}}%
        \varindent@\demoheadfont@\ignorespaces##1\unskip\definition@}%
        \rm \ignorespaces}\FN@\next@}


\newdimen\rosteritemwd
\newcount\rostercount@
\newif\iffirstitem@
\let\plainitem@\item
\newtoks\everypartoks@
\def\par@{\everypartoks@\expandafter{\the\everypar}\everypar{}}
\def\roster{\edef\leftskip@{\leftskip\the\leftskip}%
 \relaxnext@
 \rostercount@\z@  
 \def\item{\FN@\rosteritem@}%
 \DN@{\ifx\next\runinitem\let\next@\nextii@\else
  \let\next@\nextiii@\fi\next@}%
 \DNii@\runinitem  
  {\unskip  
   \DN@{\ifx\next[\let\next@\nextii@\else
    \ifx\next"\let\next@\nextiii@\else\let\next@\nextiv@\fi\fi\next@}%
   \DNii@[####1]{\rostercount@####1\relax
    \enspace{\rm(\number\rostercount@)}~\ignorespaces}%
   \def\nextiii@"####1"{\enspace{\rm####1}~\ignorespaces}%
   \def\nextiv@{\enspace{\rm(1)}\rostercount@\@ne~}%
   \par@\firstitem@false  
   \FN@\next@}%
 \def\nextiii@{\par\par@  
  \penalty\@m\smallskip\vskip-\parskip
  \firstitem@true}%
 \FN@\next@}
\def\rosteritem@{\iffirstitem@\firstitem@false\else\par\vskip-\parskip\fi
 \leftskip3\parindent\noindent  
 \DNii@[##1]{\rostercount@##1\relax
  \llap{\hbox to2.5\parindent{\hss\rm(\number\rostercount@)}%
   \hskip.5\parindent}\ignorespaces}%
 \def\nextiii@"##1"{%
  \llap{\hbox to2.5\parindent{\hss\rm##1}\hskip.5\parindent}\ignorespaces}%
 \def\nextiv@{\advance\rostercount@\@ne
  \llap{\hbox to2.5\parindent{\hss\rm(\number\rostercount@)}%
   \hskip.5\parindent}}%
 \ifx\next[\let\next@\nextii@\else\ifx\next"\let\next@\nextiii@\else
  \let\next@\nextiv@\fi\fi\next@}

\newif\ifnextRunin@
\def\endroster{\relaxnext@
 \par\leftskip@  
 \penalty-50 \vskip-\parskip\smallskip  
 \DN@{\ifx\next\Runinitem\let\next@\relax
  \else\nextRunin@false\let\item\plainitem@  
   \ifx\next\par 
    \DN@\par{\everypar\expandafter{\the\everypartoks@}}%
   \else  
    \DN@{\noindent\everypar\expandafter{\the\everypartoks@}}%
  \fi\fi\next@}%
 \FN@\next@}
\newcount\rosterhangafter@
\def\Runinitem#1\roster\runinitem{\relaxnext@
 \rostercount@\z@ 
 \def\item{\FN@\rosteritem@}%
 \def\runinitem@{#1}%
 \DN@{\ifx\next[\let\next\nextii@\else\ifx\next"\let\next\nextiii@
  \else\let\next\nextiv@\fi\fi\next}%
 \DNii@[##1]{\rostercount@##1\relax
  \def\item@{{\rm(\number\rostercount@)}}\nextv@}%
 \def\nextiii@"##1"{\def\item@{{\rm##1}}\nextv@}%
 \def\nextiv@{\advance\rostercount@\@ne
  \def\item@{{\rm(\number\rostercount@)}}\nextv@}%
 \def\nextv@{\setbox\z@\vbox  
  {\ifnextRunin@\noindent\fi  
  \runinitem@\unskip\enspace\item@~\par  
  \global\rosterhangafter@\prevgraf}%
  \firstitem@false  
  \ifnextRunin@\else\par\fi  
  \hangafter\rosterhangafter@\hangindent3\parindent
  \ifnextRunin@\noindent\fi  
  \runinitem@\unskip\enspace 
  \item@~\ifnextRunin@\else\par@\fi  
  \nextRunin@true\ignorespaces}%
 \FN@\next@}
\def\footmarkform@#1{$\m@th^{#1}$}
\let\thefootnotemark\footmarkform@
\def\makefootnote@#1#2{\insert\footins
 {\interlinepenalty\interfootnotelinepenalty
 \eightpoint\splittopskip\ht\strutbox\splitmaxdepth\dp\strutbox
 \floatingpenalty\@MM\leftskip\z@\rightskip\z@\spaceskip\z@\xspaceskip\z@
 \leavevmode{#1}\footstrut\ignorespaces#2\unskip\lower\dp\strutbox
 \vbox to\dp\strutbox{}}}
\newcount\footmarkcount@
\footmarkcount@\z@
\def\footnotemark{\let\@sf\empty\relaxnext@
 \ifhmode\edef\@sf{\spacefactor\the\spacefactor}\/\fi
 \DN@{\ifx[\next\let\next@\nextii@\else
  \ifx"\next\let\next@\nextiii@\else
  \let\next@\nextiv@\fi\fi\next@}%
 \DNii@[##1]{\footmarkform@{##1}\@sf}%
 \def\nextiii@"##1"{{##1}\@sf}%
 \def\nextiv@{\iffirstchoice@\global\advance\footmarkcount@\@ne\fi
  \footmarkform@{\number\footmarkcount@}\@sf}%
 \FN@\next@}
\def\footnotetext{\relaxnext@
 \DN@{\ifx[\next\let\next@\nextii@\else
  \ifx"\next\let\next@\nextiii@\else
  \let\next@\nextiv@\fi\fi\next@}%
 \DNii@[##1]##2{\makefootnote@{\footmarkform@{##1}}{##2}}%
 \def\nextiii@"##1"##2{\makefootnote@{##1}{##2}}%
 \def\nextiv@##1{\makefootnote@{\footmarkform@{\number\footmarkcount@}}{##1}}%
 \FN@\next@}
\def\footnote{\let\@sf\empty\relaxnext@
 \ifhmode\edef\@sf{\spacefactor\the\spacefactor}\/\fi
 \DN@{\ifx[\next\let\next@\nextii@\else
  \ifx"\next\let\next@\nextiii@\else
  \let\next@\nextiv@\fi\fi\next@}%
 \DNii@[##1]##2{\footnotemark[##1]\footnotetext[##1]{##2}}%
 \def\nextiii@"##1"##2{\footnotemark"##1"\footnotetext"##1"{##2}}%
 \def\nextiv@##1{\footnotemark\footnotetext{##1}}%
 \FN@\next@}
\def\adjustfootnotemark#1{\advance\footmarkcount@#1\relax}
\def\footnoterule{\kern-3\p@
  \hrule width 5pc\kern 2.6\p@} 
\def\captionfont@{\smc}
\def\topcaption#1#2\endcaption{%
  {\dimen@\hsize \advance\dimen@-\captionwidth@
   \rm\raggedcenter@ \advance\leftskip.5\dimen@ \rightskip\leftskip
  {\captionfont@#1}%
  \if\notempty{#2}.\enspace\ignorespaces#2\fi
  \endgraf}\nobreak\bigskip}
\def\botcaption#1#2\endcaption{%
  \nobreak\bigskip
  \setboxz@h{\captionfont@#1\if\notempty{#2}.\enspace\rm#2\fi}%
  {\dimen@\hsize \advance\dimen@-\captionwidth@
   \leftskip.5\dimen@ \rightskip\leftskip
   \noindent \ifdim\wdz@>\captionwidth@ 
   \else\hfil\fi 
  {\captionfont@#1}\if\notempty{#2}.\enspace\rm#2\fi\endgraf}}
\def\@ins{\par\begingroup\def\vspace##1{\vskip##1\relax}%
  \def\captionwidth##1{\captionwidth@##1\relax}%
  \setbox\z@\vbox\bgroup} 
\def\block{\RIfMIfI@\nondmatherr@\block\fi
       \else\ifvmode\vskip\abovedisplayskip\noindent\fi
        $$\def\endblock{\par\egroup$$}\fi
  \vbox\bgroup\advance\hsize-2\indenti\noindent}
\def\endblock{\par\egroup}
\def\cite#1{{\rm[{\citefont@\m@th#1}]}}
\def\citefont@{\rm}
\def\refsfont@{\eightpoint}
\outer\def\Refs{\runaway@{proclaim}%
 \relaxnext@ \DN@{\ifx\next\nofrills\DN@\nofrills{\nextii@}\else
  \DN@{\nextii@{References}}\fi\next@}%
 \DNii@##1{\penaltyandskip@{-200}\aboveheadskip
  \line{\hfil\headfont@\ignorespaces##1\unskip\hfil}\nobreak
  \vskip\belowheadskip
  \begingroup\refsfont@\sfcode`.=\@m}%
 \FN@\next@}
\def\endRefs{\par\endgroup}
\newbox\nobox@            \newbox\keybox@           \newbox\bybox@
\newbox\paperbox@         \newbox\paperinfobox@     \newbox\jourbox@
\newbox\volbox@           \newbox\issuebox@         \newbox\yrbox@
\newbox\pagesbox@         \newbox\bookbox@          \newbox\bookinfobox@
\newbox\publbox@          \newbox\publaddrbox@      \newbox\finalinfobox@
\newbox\edsbox@           \newbox\langbox@
\newif\iffirstref@        \newif\iflastref@
\newif\ifprevjour@        \newif\ifbook@            \newif\ifprevinbook@
\newif\ifquotes@          \newif\ifbookquotes@      \newif\ifpaperquotes@
\newdimen\bysamerulewd@
\setboxz@h{\refsfont@\kern3em}
\bysamerulewd@\wdz@
\newdimen\refindentwd
\setboxz@h{\refsfont@ 00. }
\refindentwd\wdz@
\outer\def\ref{\begingroup \noindent\hangindent\refindentwd
 \firstref@true \def\nofrills{\def\refkern@{\kern3sp}}%
 \ref@}
\def\ref@{\book@false \bgroup\let\endrefitem@\egroup \ignorespaces}
\def\moreref{\endrefitem@\endref@\firstref@false\ref@}%
\def\transl{\endrefitem@\endref@\firstref@false
  \book@false
  \prepunct@
  \setboxz@h\bgroup \aftergroup\unhbox\aftergroup\z@
    \def\endrefitem@{\unskip\refkern@\egroup}\ignorespaces}%
\def\emptyifempty@{\dimen@\wd\currbox@
  \advance\dimen@-\wd\z@ \advance\dimen@-.1\p@
  \ifdim\dimen@<\z@ \setbox\currbox@\copy\voidb@x \fi}
\let\refkern@\relax
\def\endrefitem@{\unskip\refkern@\egroup
  \setboxz@h{\refkern@}\emptyifempty@}\ignorespaces
\def\refdef@#1#2#3{\edef\next@{\noexpand\endrefitem@
  \let\noexpand\currbox@\csname\expandafter\eat@\string#1box@\endcsname
    \noexpand\setbox\noexpand\currbox@\hbox\bgroup}%
  \toks@\expandafter{\next@}%
  \if\notempty{#2#3}\toks@\expandafter{\the\toks@
  \def\endrefitem@{\unskip#3\refkern@\egroup
  \setboxz@h{#2#3\refkern@}\emptyifempty@}#2}\fi
  \toks@\expandafter{\the\toks@\ignorespaces}%
  \edef#1{\the\toks@}}
\refdef@\no{}{. }
\refdef@\key{[\m@th}{] }
\refdef@\by{}{}
\def\bysame{\by\hbox to\bysamerulewd@{\hrulefill}\thinspace
   \kern0sp}
\def\manyby{\message{\string\manyby is no longer necessary; \string\by
  can be used instead, starting with version 2.0 of \styname.STY}\by}
\refdef@\paper{\ifpaperquotes@``\fi\it}{}
\refdef@\paperinfo{}{}
\def\jour{\endrefitem@\let\currbox@\jourbox@
  \setbox\currbox@\hbox\bgroup
  \def\endrefitem@{\unskip\refkern@\egroup
    \setboxz@h{\refkern@}\emptyifempty@
    \ifvoid\jourbox@\else\prevjour@true\fi}%
\ignorespaces}
\refdef@\vol{\ifbook@\else\bf\fi}{}
\refdef@\issue{no. }{}
\refdef@\yr{}{}
\refdef@\pages{}{}
\def\page{\endrefitem@\def\pp@{\def\pp@{pp.~}p.~}\let\currbox@\pagesbox@
  \setbox\currbox@\hbox\bgroup\ignorespaces}
\def\pp@{pp.~}
\def\book{\endrefitem@ \let\currbox@\bookbox@
 \setbox\currbox@\hbox\bgroup\def\endrefitem@{\unskip\refkern@\egroup
  \setboxz@h{\ifbookquotes@``\fi}\emptyifempty@
  \ifvoid\bookbox@\else\book@true\fi}%
  \ifbookquotes@``\fi\it\ignorespaces}
\def\inbook{\endrefitem@
  \let\currbox@\bookbox@\setbox\currbox@\hbox\bgroup
  \def\endrefitem@{\unskip\refkern@\egroup
  \setboxz@h{\ifbookquotes@``\fi}\emptyifempty@
  \ifvoid\bookbox@\else\book@true\previnbook@true\fi}%
  \ifbookquotes@``\fi\ignorespaces}
\refdef@\eds{(}{, eds.)}
\def\ed{\endrefitem@\let\currbox@\edsbox@
 \setbox\currbox@\hbox\bgroup
 \def\endrefitem@{\unskip, ed.)\refkern@\egroup
  \setboxz@h{(, ed.)}\emptyifempty@}(\ignorespaces}
\refdef@\bookinfo{}{}
\refdef@\publ{}{}
\refdef@\publaddr{}{}
\refdef@\finalinfo{}{}
\refdef@\lang{(}{)}

\let\refdef@\relax 
\def\ppunbox@#1{\ifvoid#1\else\prepunct@\unhbox#1\fi}
\def\nocomma@#1{\ifvoid#1\else\changepunct@3\prepunct@\unhbox#1\fi}
\def\changepunct@#1{\ifnum\lastkern<3 \unkern\kern#1sp\fi}
\def\prepunct@{\count@\lastkern\unkern
  \ifnum\lastpenalty=0
    \let\penalty@\relax
  \else
    \edef\penalty@{\penalty\the\lastpenalty\relax}%
  \fi
  \unpenalty
  \let\refspace@\ \ifcase\count@,
\or;\or.\or 
  \or\let\refspace@\relax
  \else,\fi
  \ifquotes@''\quotes@false\fi \penalty@ \refspace@
}
\def\transferpenalty@#1{\dimen@\lastkern\unkern
  \ifnum\lastpenalty=0\unpenalty\let\penalty@\relax
  \else\edef\penalty@{\penalty\the\lastpenalty\relax}\unpenalty\fi
  #1\penalty@\kern\dimen@}
\def\endref{\endrefitem@\lastref@true\endref@
  \par\endgroup \prevjour@false \previnbook@false }
\def\endref@{%
\iffirstref@
  \ifvoid\nobox@\ifvoid\keybox@\indent\fi
  \else\hbox to\refindentwd{\hss\unhbox\nobox@}\fi
  \ifvoid\keybox@
  \else\ifdim\wd\keybox@>\refindentwd
         \box\keybox@
       \else\hbox to\refindentwd{\unhbox\keybox@\hfil}\fi\fi
  \kern4sp\ppunbox@\bybox@
\fi 
  \ifvoid\paperbox@
  \else\prepunct@\unhbox\paperbox@
    \ifpaperquotes@\quotes@true\fi\fi
  \ppunbox@\paperinfobox@
  \ifvoid\jourbox@
    \ifprevjour@ \nocomma@\volbox@
      \nocomma@\issuebox@
      \ifvoid\yrbox@\else\changepunct@3\prepunct@(\unhbox\yrbox@
        \transferpenalty@)\fi
      \ppunbox@\pagesbox@
    \fi 
  \else \prepunct@\unhbox\jourbox@
    \nocomma@\volbox@
    \nocomma@\issuebox@
    \ifvoid\yrbox@\else\changepunct@3\prepunct@(\unhbox\yrbox@
      \transferpenalty@)\fi
    \ppunbox@\pagesbox@
  \fi 
  \ifbook@\prepunct@\unhbox\bookbox@ \ifbookquotes@\quotes@true\fi \fi
  \nocomma@\edsbox@
  \ppunbox@\bookinfobox@
  \ifbook@\ifvoid\volbox@\else\prepunct@ vol.~\unhbox\volbox@
  \fi\fi
  \ppunbox@\publbox@ \ppunbox@\publaddrbox@
  \ifbook@ \ppunbox@\yrbox@
    \ifvoid\pagesbox@
    \else\prepunct@\pp@\unhbox\pagesbox@\fi
  \else
    \ifprevinbook@ \ppunbox@\yrbox@
      \ifvoid\pagesbox@\else\prepunct@\pp@\unhbox\pagesbox@\fi
    \fi \fi
  \ppunbox@\finalinfobox@
  \iflastref@
    \ifvoid\langbox@.\ifquotes@''\fi
    \else\changepunct@2\prepunct@\unhbox\langbox@\fi
  \else
    \ifvoid\langbox@\changepunct@1%
    \else\changepunct@3\prepunct@\unhbox\langbox@
      \changepunct@1\fi
  \fi
}
\outer\def\enddocument{%
 \runaway@{proclaim}%
\ifmonograph@ 
\else
 \nobreak
 \thetranslator@
 \count@\z@ \loop\ifnum\count@<\addresscount@\advance\count@\@ne
 \csname address\number\count@\endcsname
 \csname email\number\count@\endcsname
 \repeat
\fi
 \vfill\supereject\end}

\def\headfont@{\headfonts}
\def\proclaimheadfont@{\bf}
\def\specialheadfont@{\bf}
\def\subheadfont@{\bf}
\def\demoheadfont@{\smc}

\newif\ifThisToToc \ThisToTocfalse
\newif\iftocloaded \tocloadedfalse

\def\C@L{\noexpand\Cal}\def\B@B{\noexpand\Bbb}\def\fR@K{\noexpand\frak}
\def\S@{\noexpand\S}\def\P@P{\noexpand\"}
\def\xpar{\\}

\def\writetoc#1{\iftocloaded\ifThisToToc\begingroup\def\totoc{}
  \def\Cal{\noexpand\C@L}\def\Bbb{\noexpand\B@B}
  \def\frak{\noexpand\fR@K}\def\goth{\frak}\def\S{\noexpand\S@}
  \def\"{\noexpand\P@P}
  \def\xpar{\par\penalty100000 }\def\idx##1{##1}\def\\{\xpar}
  \edef\next@{\write\toc{\noindent#1\leaderfill\noexpand\folio\par}}%
  \next@\endgroup\global\ThisToTocfalse\fi\fi}
\def\leaderfill{\leaders\hbox to 1em{\hss.\hss}\hfill}

\newif\ifindexloaded \indexloadedfalse
\def\idx#1{\ifindexloaded\begingroup\def\ign{}\def\it{}\def\/{}%
 \def\smc{}\def\bf{}\def\tt{}%
 \def\Cal{\noexpand\C@L}\def\Bbb{\noexpand\B@B}%
 \def\frak{\noexpand\fR@K}\def\goth{\frak}\def\S{\noexpand\S@}%
  \def\"{\noexpand\P@P}%
 {\edef\next@{\write\index{#1, \noexpand\folio}}\next@}%
 \endgroup\fi{#1}}
\def\ign#1{}

\def\totoc{\global\ThisToToctrue}

\outer\def\head#1\endhead{\par\penaltyandskip@{-200}\aboveheadskip
 {\headfont@\raggedcenter@\interlinepenalty\@M
 \ignorespaces#1\endgraf}\nobreak
 \vskip\belowheadskip
 \headmark{#1}\writetoc{#1}}

\outer\def\chaphead#1\endchaphead{\par\penaltyandskip@{-200}\aboveheadskip
 {\chapheadfonts\raggedcenter@\interlinepenalty\@M
 \ignorespaces#1\endgraf}\nobreak
 \vskip3\belowheadskip
 \headmark{#1}\writetoc{#1}}

\def\folio{{\foliofont@\ifnum\pageno<\z@ \romannumeral-\pageno
 \else\number\pageno \fi}}
\newtoks\leftheadtoks
\newtoks\rightheadtoks

\def\leftheadtext{\nofrills@{\relax}\lht@
  \DNii@##1{\leftheadtoks\expandafter{\lht@{##1}}%
    \mark{\the\leftheadtoks\noexpand\else\the\rightheadtoks}
    \ifsyntax@\setboxz@h{\def\\{\unskip\space\ignorespaces}%
        \headlinefont@##1}\fi}%
  \FN@\next@}
\def\rightheadtext{\nofrills@{\relax}\rht@
  \DNii@##1{\rightheadtoks\expandafter{\rht@{##1}}%
    \mark{\the\leftheadtoks\noexpand\else\the\rightheadtoks}%
    \ifsyntax@\setboxz@h{\def\\{\unskip\space\ignorespaces}%
        \headlinefont@##1}\fi}%
  \FN@\next@}
\def\NoRunningHeads{\global\runheads@false\global\let\headmark\eat@}

\newif\iffirstpage@     \firstpage@true
\newif\ifrunheads@      \runheads@true

\newdimen\fullhsize \fullhsize=\hsize
\newdimen\fullvsize \fullvsize=\vsize
\def\fullline{\hbox to\fullhsize}

\def\pagenumbers{\gdef\folio{\folio@}}

\let\norunningheads\NoRunningHeads
\def\userunningheads{\global\runheads@true}
\norunningheads

\headline={\def\chapter#1{\chapterno@. }%
  \def\\{\unskip\space\ignorespaces}\ifrunheads@\headlinefont@
    \ifodd\pageno\rightheadline \else\leftheadline\fi
   \else\hfil\fi\ifNoRunHeadline\global\NoRunHeadlinefalse\fi}
\let\folio@\folio
\def\foliofont@{\foliofont}
\def\foliofont{\eightrm}
\def\headlinefont@{\headlinefont}
\def\headlinefont{\eightpoint\smc}
\def\leftheadline{\rlap{\folio}\hfill
   \ifNoRunHeadline\else\iftrue\topmark\fi\fi \hfill}
\def\rightheadline{\hfill\ifNoRunHeadline
   \else \expandafter\fi
  \hfill \llap{\folio}}
\footline={{\eightpoint\bottremark}%
   \ifrunheads@\else\hfil{\let\foliofont\tenrm\folio}\fi\hfil}
\def\bottremark{}
 
\newif\ifNoRunHeadline      
\def\norunninghead{\global\NoRunHeadlinetrue}
\norunninghead

\output={\output@}
%
\newif\ifoffset\offsetfalse
\output={\output@}
\def\output@{%
 \ifoffset 
  \ifodd\count0\advance\hoffset by0.5truecm
   \else\advance\hoffset by-0.5truecm\fi\fi
 \shipout\vbox{%
  \makeheadline \pagebody \makefootline }%
 \advancepageno \ifnum\outputpenalty>-\@MM\else\dosupereject\fi}

\def\indexoutput#1{%
  \ifoffset 
   \ifodd\count0\advance\hoffset by0.5truecm
    \else\advance\hoffset by-0.5truecm\fi\fi
  \shipout\vbox{\makeheadline
  \vbox to\fullvsize{\boxmaxdepth\maxdepth%
  \ifvoid\topins\else\unvbox\topins\fi%
  #1 %
  \ifvoid\footins\else 
    \vskip\skip\footins
    \footnoterule
    \unvbox\footins\fi
  \ifr@ggedbottom \kern-\dimen@ \vfil \fi}%
  \baselineskip2pc
  \makefootline}%
 \global\advance\pageno\@ne
 \ifnum\outputpenalty>-\@MM\else\dosupereject\fi}
 
 \newbox\partialpage \newdimen\halfsize \halfsize=0.5\fullhsize
 \advance\halfsize by-0.5em

 \def\begindoublecolumns{\output={\indexoutput{\unvbox255}}%
   \begingroup \def\line{\fullline}
   \output={\global\setbox\partialpage=\vbox{\unvbox255\bigskip}}\eject
   \output={\doublecolumnout}\hsize=\halfsize \vsize=2\fullvsize}
 \def\enddoublecolumns{\output={\balancecolumns}\eject
  \endgroup \pagegoal=\fullvsize%
  \output={\output@}}
\def\doublecolumnout{\splittopskip=\topskip \splitmaxdepth=\maxdepth
  \dimen@=\fullvsize \advance\dimen@ by-\ht\partialpage
  \setbox0=\vsplit255 to \dimen@ \setbox2=\vsplit255 to \dimen@
  \indexoutput{\pagesofar} \unvbox255 \penalty\outputpenalty}
\def\pagesofar{\unvbox\partialpage
  \wd0=\hsize \wd2=\hsize \hbox to\fullhsize{\box0\hfil\box2}}
\def\balancecolumns{\setbox0=\vbox{\unvbox255} \dimen@=\ht0
  \advance\dimen@ by\topskip \advance\dimen@ by-\baselineskip
  \divide\dimen@ by2 \splittopskip=\topskip
  {\vbadness=10000 \loop \global\setbox3=\copy0
    \global\setbox1=\vsplit3 to\dimen@
    \ifdim\ht3>\dimen@ \global\advance\dimen@ by1pt \repeat}
  \setbox0=\vbox to\dimen@{\unvbox1} \setbox2=\vbox to\dimen@{\unvbox3}
  \pagesofar}

\tenpoint
\catcode`\@=\active

\def\smallheadings{\let\chapheadfonts\tenpoint\let\headfonts\tenpoint}

\tenpoint
\catcode`\@=\active

\def\LL{\leavevmode\setbox0=\hbox{L}\hbox to\wd0{\hss\char'40L}}
\def\al{\alpha}
\def\be{\beta}

\def\de{\delta}
\def\ep{\varepsilon}

\def\si{\sigma}

\def\De{\Delta}


\def\today{\ifcase\month\or
 January\or February\or March\or April\or May\or June\or
 July\or August\or September\or October\or November\or December\fi
 \space\number\day, \number\year}

\def\({\left(}
\def\){\right)}
\def\[{\left[}
\def\]{\right]}

\def\3{\ss}
\catcode`\@=11
\def\dddot#1{\vbox{\ialign{##\crcr
      .\hskip-.5pt.\hskip-.5pt.\crcr\noalign{\kern1.5\p@\nointerlineskip}
      $\hfil\displaystyle{#1}\hfil$\crcr}}}

\newif\iftab@\tab@false
\newif\ifvtab@\vtab@false
\def\tab{\bgroup\tab@true\vtab@false\vst@bfalse\Strich@false%
   \def\\{\global\hline@@false%
     \ifhline@\global\hline@false\global\hline@@true\fi\cr}
   \edef\l@{\the\leftskip}\ialign\bgroup\hskip\l@##\hfil&&##\hfil\cr}
\def\endtab{\cr\egroup\egroup}
\def\vtab{\vtop\bgroup\vst@bfalse\vtab@true\tab@true\Strich@false%
   \bgroup\def\\{\cr}\ialign\bgroup&##\hfil\cr}
\def\endvtab{\cr\egroup\egroup\egroup}
\def\stab{\D@cke0.5pt\null 
 \bgroup\tab@true\vtab@false\vst@bfalse\Strich@true\Let@@\vspace@
 \normalbaselines\offinterlineskip
  \openup\spreadmlines@
 \edef\l@{\the\leftskip}\ialign
 \bgroup\hskip\l@##\hfil&&##\hfil\crcr}
\def\endstab{\crcr\egroup
 \egroup}
\newif\ifvst@b\vst@bfalse
\def\vstab{\D@cke0.5pt\null
 \vtop\bgroup\tab@true\vtab@false\vst@btrue\Strich@true\bgroup\Let@@\vspace@
 \normalbaselines\offinterlineskip
  \openup\spreadmlines@\bgroup}
\def\endvstab{\crcr\egroup\egroup
 \egroup\tab@false\Strich@false}

\newdimen\htstrut@
\htstrut@8.5\p@
\newdimen\htStrut@
\htStrut@12\p@
\newdimen\dpstrut@
\dpstrut@3.5\p@
\newdimen\dpStrut@
\dpStrut@3.5\p@
\def\openup{\afterassignment\@penup\dimen@=}
\def\@penup{\advance\lineskip\dimen@
  \advance\baselineskip\dimen@
  \advance\lineskiplimit\dimen@
  \divide\dimen@ by2
  \advance\htstrut@\dimen@
  \advance\htStrut@\dimen@
  \advance\dpstrut@\dimen@
  \advance\dpStrut@\dimen@}
\def\Let@@{\relax%
    \def\\{\global\hline@@false%
     \ifhline@\global\hline@false\global\hline@@true\fi\cr}%
    \iffalse}\fi}
\def\matrix{\null\,\vcenter\bgroup
 \tab@false\vtab@false\vst@bfalse\Strich@false\Let@@\vspace@
 \normalbaselines\openup\spreadmlines@\ialign
 \bgroup\hfil$\m@th##$\hfil&&\quad\hfil$\m@th##$\hfil\crcr
 \Mathstrut@\crcr\noalign{\kern-\baselineskip}}
\def\endmatrix{\crcr\Mathstrut@\crcr\noalign{\kern-\baselineskip}\egroup
 \egroup\,}
\def\smatrix{\D@cke0.5pt\null\,
 \vcenter\bgroup\tab@false\vtab@false\vst@bfalse\Strich@true\Let@@\vspace@
 \normalbaselines\offinterlineskip
  \openup\spreadmlines@\ialign
 \bgroup\hfil$\m@th##$\hfil&&\quad\hfil$\m@th##$\hfil\crcr}
\def\endsmatrix{\crcr\egroup
 \egroup\,\Strich@false}
\newdimen\D@cke
\def\Dicke#1{\global\D@cke#1}
\newtoks\tabs@\tabs@{&}
\newif\ifStrich@\Strich@false
\newif\iff@rst

\def\Stricherr@{\iftab@\ifvtab@\errmessage{\noexpand\s not allowed
     here. Use \noexpand\vstab!}%
  \else\errmessage{\noexpand\s not allowed here. Use \noexpand\stab!}%
  \fi\else\errmessage{\noexpand\s not allowed
     here. Use \noexpand\smatrix!}\fi}
\def\format{\ifvst@b\else\crcr\fi\egroup\iffalse{\fi\ifnum`}=0 \fi\format@}
\def\format@#1\\{\def\preamble@{#1}%
 \def\Str@chfehlt##1{\ifx##1\s\Stricherr@\fi\ifx##1\\\let\Next\relax%
   \else\let\Next\Str@chfehlt\fi\Next}%
 \def\c{\hfil\noexpand\ifhline@@\hbox{\vrule height\htStrut@%
   depth\dpstrut@ width\z@}\noexpand\fi%
   \ifStrich@\hbox{\vrule height\htstrut@ depth\dpstrut@ width\z@}%
   \fi\iftab@\else$\m@th\fi\the\hashtoks@\iftab@\else$\fi\hfil}%
 \def\r{\hfil\noexpand\ifhline@@\hbox{\vrule height\htStrut@%
   depth\dpstrut@ width\z@}\noexpand\fi%
   \ifStrich@\hbox{\vrule height\htstrut@ depth\dpstrut@ width\z@}%
   \fi\iftab@\else$\m@th\fi\the\hashtoks@\iftab@\else$\fi}%
 \def\l{\noexpand\ifhline@@\hbox{\vrule height\htStrut@%
   depth\dpstrut@ width\z@}\noexpand\fi%
   \ifStrich@\hbox{\vrule height\htstrut@ depth\dpstrut@ width\z@}%
   \fi\iftab@\else$\m@th\fi\the\hashtoks@\iftab@\else$\fi\hfil}%
 \def\s{\ifStrich@\ \the\tabs@\vrule width\D@cke\the\hashtoks@%
          \fi\the\tabs@\ }%
 \def\sa{\ifStrich@\vrule width\D@cke\the\hashtoks@%
            \the\tabs@\ %
            \fi}%
 \def\se{\ifStrich@\ \the\tabs@\vrule width\D@cke\the\hashtoks@\fi}%
 \def\cd{\hfil\noexpand\ifhline@@\hbox{\vrule height\htStrut@%
   depth\dpstrut@ width\z@}\noexpand\fi%
   \ifStrich@\hbox{\vrule height\htstrut@ depth\dpstrut@ width\z@}%
   \fi$\dsize\m@th\the\hashtoks@$\hfil}%
 \def\rd{\hfil\noexpand\ifhline@@\hbox{\vrule height\htStrut@%
   depth\dpstrut@ width\z@}\noexpand\fi%
   \ifStrich@\hbox{\vrule height\htstrut@ depth\dpstrut@ width\z@}%
   \fi$\dsize\m@th\the\hashtoks@$}%
 \def\ld{\noexpand\ifhline@@\hbox{\vrule height\htStrut@%
   depth\dpstrut@ width\z@}\noexpand\fi%
   \ifStrich@\hbox{\vrule height\htstrut@ depth\dpstrut@ width\z@}%
   \fi$\dsize\m@th\the\hashtoks@$\hfil}%
 \ifStrich@\else\Str@chfehlt#1\\\fi%
 \setbox\z@\hbox{\xdef\Preamble@{\preamble@}}\ifnum`{=0 \fi\iffalse}\fi
 \ialign\bgroup\span\Preamble@\crcr}
\newif\ifhline@\hline@false
\newif\ifhline@@\hline@@false
\def\hlinefor#1{\multispan@{\strip@#1 }\leaders\hrule height\D@cke\hfill%
    \global\hline@true\ignorespaces}
\def\Item "#1"{\par\noindent\hangindent2\parindent%
  \hangafter1\setbox0\hbox{\rm#1\enspace}\ifdim\wd0>2\parindent%
  \box0\else\hbox to 2\parindent{\rm#1\hfil}\fi\ignorespaces}
\def\ITEM #1"#2"{\par\noindent\hangafter1\hangindent#1%
  \setbox0\hbox{\rm#2\enspace}\ifdim\wd0>#1%
  \box0\else\hbox to 0pt{\rm#2\hss}\hskip#1\fi\ignorespaces}
\def\item"#1"{\par\noindent\hang%
  \setbox0=\hbox{\rm#1\enspace}\ifdim\wd0>\the\parindent%
  \box0\else\hbox to \parindent{\rm#1\hfil}\enspace\fi\ignorespaces}
\let\plainitem@\item
\catcode`\@=13

\hsize13cm
\vsize19cm

\magnification1200

\TagsOnRight

\def\ArmDAA{1}
\def\AthaAI{2}
\def\AtReAA{3}
\def\BesDAA{4}
\def\BRWaAA{5}
\def\BRWaAB{6}
\def\CarlAP{7}
\def\ChaFAA{8}
\def\EdelAA{9}
\def\FoReAA{10}
\def\FSZeAC{11}
\def\FSZeAD{12}
\def\FOZeAB{13}
\def\GrKPAA{14}
\def\HumpAC{15}
\def\KrewAC{16}
\def\ReivAG{17}
\def\SlatAC{18}
\def\StanAP{19}
\def\StemAZ{20}
\def\TzanAA{21}

\def\aa{1}
\def\AA{2}
\def\AB{3}
\def\AC{4}
\def\AD{5}
\def\BDa{6}
\def\BD{7}
\def\BH{8}
\def\BA{9}
\def\BB{10}
\def\BC{11}
\def\BE{12}
\def\BF{13}
\def\BG{14}
\def\BI{15}
\def\BJ{16}
\def\BK{17}
\def\BKa{17a}
\def\BKb{17b}
\def\BKc{17c}
\def\BKd{17d}
\def\BL{18}
\def\BM{19}
\def\BMa{20}
\def\BN{21}
\def\CA{22}
\def\CB{23}
\def\CC{24}
\def\CD{25}
\def\CE{26}
\def\CF{27}
\def\CG{28}
\def\CH{29}
\def\BO{30}
\def\BP{31}
\def\BPa{32}
\def\BQ{33}
\def\BQa{34}
\def\BR{35}
\def\BS{36}
\def\BSa{36a}
\def\BSb{36b}
\def\BSc{36c}
\def\BT{37}
\def\BU{38}
\def\BV{39}
\def\BX{40}
\def\AE{41}
\def\aE{42}
\def\Ae{43}
\def\Aee{44}
\def\Aa{45}
\def\Ab{46}
\def\Ac{47}
\def\Aeee{48}
\def\BY{49}
\def\Af{50}
\def\AF{51}
\def\AG{52}
\def\AH{53}
\def\BZ{54}

\def\al{\alpha}
\def\be{\beta}

\def\rk{\operatorname{rk}}
\def\Nar{\operatorname{Nar}}

\topmatter 
\title The $F$-triangle of the generalised cluster complex
\endtitle 
\author C.~Krattenthaler$^{\dagger\ddagger}$
\endauthor 
\affil 
Institut Camille Jordan, Universit\'e Claude Bernard Lyon-I,\\
21, avenue Claude Bernard, F-69622 Villeurbanne Cedex, France.\\
WWW: \tt http://igd.univ-lyon1.fr/\~{}kratt
\endaffil
\address Institut Camille Jordan, Universit\'e Claude Bernard Lyon-I, 
21, avenue Claude Bernard, F-69622 Villeurbanne Cedex, France.
WWW: \tt http://igd.univ-lyon1.fr/\~{}kratt
\endaddress
\thanks{$^\dagger$Research partially supported by
EC's IHRP Programme,
grant HPRN-CT-2001-00272, ``Algebraic Combinatorics in Europe."\newline
\indent$^\ddagger$Current address: Fakult\"at f\"ur Mathematik, 
Universit\"at Wien, Nordbergstra{\ss}e~15, A-1090 Vienna, Austria}%
\endthanks
\subjclass Primary 05E15;
 Secondary 05A05 05A10 05A15 05A19 06A07 20F55 33C05 
\endsubjclass
\keywords root systems, 
generalised cluster complex, face numbers, $F$-triangle,
$m$-divisible non-crossing partitions, M\"obius function,
$M$-triangle, Chu--Vandermonde summation.
\endkeywords
\abstract 
The $F$-triangle is a refined face count for the generalised cluster
complex of Fomin and
Reading. We compute the $F$-triangle explicitly for all irreducible
finite root systems. Furthermore, we use these results to partially
prove the ``$F=M$ Conjecture" of Armstrong which predicts a surprising
relation between the $F$-triangle and the M\"obius function of
his $m$-divisible partition poset associated to a 
finite root system.
\endabstract
\endtopmatter
\document

\subhead 1. Introduction\endsubhead
Fomin and Zelevinsky created a new exciting research field when they invented
{\it cluster algebras} in \cite{\FSZeAC}. The classification of
cluster algebras {\it of finite type} from \cite{\FSZeAD} 
says that there is a one-to-one correspondence between
finite-type cluster algebras and finite root systems. 
Furthermore, for each finite root system $\Phi$, Fomin and Zelevinsky
\cite{\FOZeAB} defined a simplicial complex corresponding to the associated
cluster algebra, the {\it cluster complex} $\De(\Phi)$. This is
a simplicial complex on a subset of the set of roots $\Phi$. As they
showed, this complex has many remarkable properties. In particular, 
the number of facets is given by the {\it Catalan number for the root
system $\Phi$}, and, moreover, all the face numbers are given by elegant
product formulae. Further remarkable (originally, conjectural)
properties have been discovered by Chapoton in \cite{\ChaFAA}.
In this paper, he refines the face enumeration to, what he calls,
the  ``$F$-triangle." He computed the $F$-triangle for all types
(and revealed his findings partially in \cite{\ChaFAA}) and
observed a surprising relationship (see \cite{\ChaFAA, Conjecture~1}) 
between the $F$-triangle and the M\"obius
function of the non-crossing partition lattice $NC(\Phi)$ associated
to $\Phi$, the latter being
due to Bessis \cite{\BesDAA} and Brady and Watt \cite{\BRWaAA}.
This relationship, to which we shall refer in the sequel as the
``$F=M$ Conjecture," has been recently proved by Athanasiadis
\cite{\AthaAI}. For further fascinating properties of the $F$-triangle
see \cite{\ChaFAA}. 

The subject of this paper is ``$m$-generalisations" of the cluster
complex $\De(\Phi)$ and of the non-crossing partition lattice
$NC(\Phi)$. More precisely, in \cite{\FoReAA}, Fomin and Reading
introduce the generalised cluster complex $\De^m(\Phi)$, where $m$ is
some non-negative integer. This is, again, a simplicial complex, now
on {\it coloured\/} roots, 
and for $m=1$ it reduces to the (ordinary) cluster complex
$\De(\Phi)$. As they show, this generalised complex has again
remarkable properties. In particular, the number of facets is given by
the {\it Fuss--Catalan number for the root
system $\Phi$}, and, moreover, again, all the face numbers are
given by elegant product formulae. 

Going one step further, Sergey Fomin suggested to the author to
investigate 
the ``Cha\-po\-ton-like" refinement of the face numbers of the generalised
cluster complex $\De^m(\Phi)$, that is to say, to study the ``$F$-triangle"
for $\De^m(\Phi)$ (see Section~2 for the definition). 
This is what we do in this paper.
We compute the $F$-triangle of $\De^m(\Phi)$ for all types of
irreducible root systems $\Phi$, see Sections~4--7. We do this
case-by-case. While a case-independent formula would certainly be
desirable, certain features of our results (in particular, the
appearance of the
Kronecker delta in the formula in Theorem~FD in Section~6 
for type $D_n$) make it highly unlikely
that such a case-independent formula exists.
As an aside, we draw the reader's attention to the unexpected outcome
of our results that the refined face numbers are all polynomials in
$m$ {\it with non-negative coefficients}, a phenomenon for which we
have no intrinsic explanation.

One may then ask if there is also an ``$F=M$ Conjecture" in this
generalised context. This would, first of all, require an
``$m$-extension" of the non-crossing partition lattice. Indeed,
Armstrong \cite{\ArmDAA} has recently introduced the
``$m$-divisible non-crossing partition poset" $NC^m(\Phi)$,
generalising an earlier construction of Edelman \cite{\EdelAA} in type
$A_n$. He
shows that this poset has also remarkable properties, resembling those
of the non-crossing partition lattices.
Moreover, he observed that there is also a rather
straight-forward extension of the $F=M$ Conjecture relating the
$F$-triangle of the generalised cluster complex $\De^m(\Phi)$
to the M\"obius function of the corresponding $m$-divisible
non-crossing partitions poset $NC^m(\Phi)$. We reproduce this conjecture in
Section~8. (We refer the reader to \cite{\ArmDAA, Sec.~4}
and \cite{\TzanAA} for further
fascinating properties of the $F$-triangle of $\De^m(\Phi)$.)

With the explicit formulae for the $F$-triangle in hand, we are able
to prove this ``$m$-version" of the $F=M$ Conjecture in types
$A_n$ and $B_n$, for the dihedral root systems $I_2(a)$, for
the hyperbolic root systems $H_3$ and $H_4$, and for $F_4$ and $E_6$, 
see Sections~9, 10, 13--17. In types $A_n$ and $B_n$, the
proofs depend crucially on results about rank selected chain
enumeration in $NC^m(A_n)$ and $NC^m(B_n)$ due to Edelman
\cite{\EdelAA} and Armstrong \cite{\ArmDAA}, respectively.
Moreover, in Section~11, we provide a calculation in type $D_n$ which will
prove the conjecture also for this type once the corresponding 
rank selected chain
enumeration result analogous to 
the ones by Edelman and Armstrong
is available for $NC^m(D_n)$. In view of the results of Athanasiadis
and Reiner \cite{\AtReAA} on rank selected
chain enumeration in $NC^1(D_n)=NC(D_n)$, this argument does
accomplish the proof for $m=1$.
As we explain in Section~12, the
verification of the (generalised) $F=M$ Conjecture  in the 
exceptional types is a routine task which can, in principle, be
carried out on a computer. To do this in practice for the root systems
$E_7$ or $E_8$, say, may however require additional simplifications
of the proposed procedure.\footnote{Using some additional ideas, this
has been carried through in {\it``The $M$-triangle of generalised non-crossing 
partitions for the types $E_7$ and $E_8$"} (preprint; 
{\tt arXiv:math.CO/0601676}).}

In the final Section~18, we prove Armstrong's conjecture
\cite{\ArmDAA, Sec.~4} on the form of, what he calls, the {\it dual
$F$-triangle} in a case-by-case fashion. For the exceptional root
systems this is just a routine calculation, while for $A_n$, $B_n$
and $D_n$ this requires only the Chu--Vandermonde summation formula.

We conclude the introduction by saying a few words how the
$F$-triangles for $\De^m(\Phi)$ are found and the corresponding
results are proved in this paper.
The main tool in \cite{\ChaFAA} for finding $F$-triangles for the
(ordinary) cluster complex $\De(\Phi)$ consists in two recurrence formulae (see
\cite{\ChaFAA, Prop.~3}). These recurrences carry over verbatim 
to the generalised cluster complex $\De^m(\Phi)$, see Proposition~F in
Section~2. 
Indeed, in the exceptional types, the formula for the $F$-triangle 
for $\De^m(\Phi)$ can be
found in a routine fashion by using these recurrences, see Section~7.
In types
$A_n$, $B_n$, and $D_n$, however, the recurrences can only be used to
compute the $F$-triangle for the corresponding generalised cluster
complex for {\it specific} $n$. By doing this for
sufficiently many $n$, we first worked out {\it guesses} for the
$F$-triangle for {\it generic} $n$. (In types $A_n$ and $B_n$, this has also
been done independently by Tzanaki \cite{\TzanAA}.)
Subsequently, one tries to verify these
guesses by checking the recurrences. As it turns out, this requires 
multivariate summation formulae due to Carlitz \cite{\CarlAP},
which we restate here in Section~3 for the convenience of the reader.
An interesting detail is the fact that Chapoton's proofs 
in \cite{\ChaFAA} for the
$F$-triangle for $\De(A_n)$ and $\De(B_n)$, respectively, which also
use Carlitz's summation formulae, do {\it not\/} extend to
$\De^m(A_n)$ and $\De^m(B_n)$, for the following reason. 
In order to do the above described verification using the recurrences,
he has to evaluate a triple sum. He does this by first simplifying
one sum by means of the Chu--Vandermonde summation, and by using
subsequently one of Carlitz's summation formulae to evaluate the
remaining double sum. However, if $m\ne1$, the Chu--Vandermonde
summation is not applicable to the triple sum that we encounter at the
start. Remarkably, it is
possible to apply Carlitz's summation formula {\it directly}, 
in a {\it different\/} way than in \cite{\ChaFAA}.
The use of the Chu--Vandermonde summation is
then not necessary anymore.

\subhead 2. Preliminaries\endsubhead
Let $\Phi$ be a finite root system of rank $n$. (We refer the reader to
\cite{\HumpAC} for all root system terminology.)
For a non-negative integer $m$, the generalised cluster complex 
$\De^m(\Phi)$ is a certain simplicial complex on a certain set of
``coloured" roots, the roots being from $\Phi$. The precise definition
will not be important here, we refer the reader to \cite{\FoReAA, Sec.~2}.
The only fact which is important here is that some of the coloured
roots can be positive, others negative. Let $f_{k,l}(\Phi,m)$ denote
the number of faces of $\De^m(\Phi)$ which contain exactly $k$
positive and $l$ negative coloured roots. Define the $F$-triangle of
$\De^m(\Phi)$, denoted by $F^m_\Phi(x,y)$, as the two-variable polynomial
$$F^m_\Phi(x,y)=
\sum _{k,l\ge0} ^{}f_{k,l}(\Phi,m)\,x^ky^l.
\tag\aa$$
It is called ``triangle" because all faces have cardinality at most
$n$ and, thus, in the summation in (\aa) we can restrict the summation
indices to the triangle $k+l\le n$, $k,l\ge0$.

Then, in this generalised
context, the arguments from \cite{\ChaFAA, Prop.~3} carry over
verbatim to prove the following properties of the $F$-triangle of
$\De^m(\Phi)$.

\proclaim{Proposition F}
The $F$-triangle $F^m_\Phi(x,y)$ satisfies the following three properties:

\roster 
\item If $\Phi$ and $\Phi'$ are two root systems, then
$$F^m_{\Phi\times\Phi'}(x,y)=F^m_\Phi(x,y)F^m_{\Phi'}(x,y),
\tag\AA$$
where $\Phi\times\Phi'$ denotes the orthogonal product of the two root
systems.
\item If $\Phi=\Phi(S)$ is an irreducible root system with simple roots $S$, then
$$
\frac {\partial} {\partial y}F^m_{\Phi(S)}(x,y)=
\sum _{\al\in S} ^{}F^m_{\Phi(S\backslash\{\al\})}(x,y),
\tag\AB$$
where $\Phi(S\backslash\{\al\})$ denotes the root system generated by
the simple roots $S\backslash\{\al\}$.
\item The specialisation $x=y$ is given by
$$F^m_{\Phi}(x,x)=
\sum _{k\ge0} ^{}f_k(\Phi,m)\,x^k,
\tag\AC$$
where the coefficients 
$f_k(\Phi,m)$ are the face numbers of the cluster complex
$\De^m(\Phi)$, summarised in \cite{\FoReAA, Theorem~7.5} for the
irreducible root systems.
\endroster
\endproclaim

We remark that an equivalent statement of (\AB) is
$$l\cdot f_{k,l}(\Phi(S),m)=
\sum _{\al\in S} ^{}
f_{k,l-1}(\Phi(S\backslash \{\al\}),m),\quad k,l\ge0.
\tag\AD$$
Moreover, in view of the multiplicativity property (\AA), 
it suffices to compute the $F$-triangle for the irreducible root systems,
which we do in Sections~4--7.

\subhead 3. Carlitz's summation formulae\endsubhead
Crucial in the proofs of our claims for the $F$-triangle in types
$A_n$, $B_n$ and $D_n$ are the following two 
double sum evaluations due to Carlitz \cite{\CarlAP}. (He has
in fact extensions for any number of summations, see \cite{\CarlAP, Sec.~6}.)
Let
$$A_{k,n}(\al,\be)=\frac {bk\al+cn\be+\al\be} {(ak+cn+\al)(bk+dn+\be)}
\binom {ak+cn+\al}k\binom {bk+dn+\be}n.$$
Here, and in the sequel, for integers $N$ and $K$
the binomial coefficient $\binom NK$ is understood
according to the definition
$$\binom NK=\cases \frac {N(N-1)\cdots(N-K+1)} {K!}&\text{if }K\ge0,\\
0&\text{if }K<0.\endcases\tag\BDa$$
Then (see \cite{\CarlAP, (5.14)}),
$$\sum _{k_1,n_1\ge0}
^{}A_{k_1,n_1}(\al,\be)A_{k-k_1,n-n_1}(\al',\be')=
A_{k,n}(\al+\al',\be+\be').\tag\BD$$
Furthermore (see 
\cite{\CarlAP, (5.15); the minus sign in front of $cn$ must be
replaced by a plus sign there}),
$$\multline 
\sum _{k_1,n_1\ge0}
^{}\binom {ak_1+cn_1+\al-1}{k_1}\binom {bk_1+dn_1+\be-1}{n_1}
A_{k-k_1,n-n_1}(\al',\be')\\=
\binom {ak+cn+\al+\al'-1}k \binom
{bk+dn+\be+\be'-1}{n}.
\endmultline\tag\BH$$

\subhead 4. The $F$-triangle for $A_n$\endsubhead
The theorem below gives an explicit expression for 
the refined face numbers $f_{k,l}(A_n,m)$, and, thus, of the
$F$-triangle in type $A_n$.

\proclaim{Theorem FA}
For $n\ge1$, the face numbers $f_{k,l}(A_n,m)$ are given by
$$f_{k,l}(A_n,m)=\frac{l+1 }{k+l+1}{\binom {n} { k+l}} 
{\binom {m (n+1)+k-1} k}.$$
\endproclaim

\demo{Proof}
In view of Proposition~F and (\AD) in Section~2,
in order to prove this claim we have to show
$$l\cdot f_{k,l}(A_n,m)=
\underset l_1+l_2=l-1\to{\underset k_1+k_2=k\hphantom{-1}
\to{\sum _{n_1+n_2=n-1} ^{}}}
f_{k_1,l_1}(A_{n_1},m)\cdot f_{k_2,l_2}(A_{n_2},m)
\tag\BA$$
and
$$
\sum _{k_1+l_1=k} ^{}f_{k_1,l_1}(A_n,m)=
\frac {1} {k+1}\binom {n}k\binom {m(n+1)+k+1}k.
\tag\BB$$

The triple sum on the right-hand side of (\BA) is 
$$\multline 
\sum _{k_1,l_1,n_1\ge0} ^{}\frac{l_1+1 }{k_1+l_1+1}{\binom
{n_1} { k_1+l_1}} 
{\binom {m (n_1+1)+k_1-1} {k_1}}\\
\cdot
\frac{l-l_1 }{k-k_1+l-l_1}{\binom {n-n_1-1} {
k-k_1+l-l_1-1}} 
{\binom {m (n-n_1)+k-k_1-1} {k-k_1}}.
\endmultline$$
We replace $n_1$ by $n_1+k_1+l_1$ and rewrite the resulting
expression in the form
$$\multline 
\sum _{k_1,l_1,n_1\ge0} ^{}\frac{m(l_1+1) }{m (k_1+l_1+n_1+1)+k_1}{\binom
{n_1+k_1+l_1+1} {n_1}}  
{\binom {m (n_1+k_1+l_1+1)+k_1} {k_1}}\\
\cdot
\frac{m(l-l_1) }{m (n-n_1-k_1-l_1)+k-k_1}
{\binom {n-n_1-k_1-l_1} {n-k-l-n_1}} \\
\cdot
{\binom {m (n-n_1-k_1-l_1)+k-k_1} {k-k_1}}.
\endmultline\tag\BC$$
Forgetting the sum over $l_1$, this is now exactly in the form
of the left-hand side of (\BD) with $n$ replaced by $n-k-l$,
$a=m+1$, $c=m$, $\al=m(l_1+1)$, $\al'=m(l-l_1)$, $b=d=1$,
$\be=l_1+1$, and $\be'=l-l_1$.
Substituting the right-hand side, we obtain
$$
\sum _{l_1=0} ^{l-1}A_{k,n-k-l}(m(l+1),l+1)$$
for the sum in (\BC), or, equivalently,
$$l\cdot \frac{m(l+1) }{m (n+1)+k}{\binom {m (n+1)+k} k}
{\binom {n+1} { n-k-l}} ,$$
which is indeed equal to $l\cdot f_{k,l}(A_n,m)$.
This proves (\BA).

\medskip
Next we compute the sum on the left-hand side of (\BB),
$$\align 
\sum _{k_1=0} ^{k}&\frac{k-k_1+1 }{k+1}{\binom {n} { k}} 
{\binom {m (n+1)+k_1-1} {k_1}}\\
&=\sum _{k_1=0} ^{k}{\binom {n} { k}} 
{\binom {m (n+1)+k_1-1} {k_1}}
-\sum _{k_1=0} ^{k}\frac{m(n+1) }{k+1}{\binom {n} { k}} 
{\binom {m (n+1)+k_1-1} {k_1-1}}\\
&={\binom {n} { k}} 
\binom {m (n+1)+k} {m (n+1)}
-\frac{m(n+1) }{k+1}{\binom {n} { k}} 
\binom {m (n+1)+k} {m (n+1)+1}\\
&=\frac {1} {k+1}{\binom {n} { k}} 
\binom {m (n+1)+k+1} {k},
\endalign$$
the simplification of summations being due to the Chu--Vandermonde
summation (see \cite{\GrKPAA, Sec.~5.1, (5.27)} or \cite{\SlatAC, (1.7.7),
Appendix~(III.4)}).
This completes the proof.\quad \quad \qed
\enddemo

\subhead 5. The $F$-triangle for $B_n$\endsubhead
The theorem below gives an explicit expression for 
the refined face numbers $f_{k,l}(B_n,m)$, and, thus, of the
$F$-triangle in type $B_n$.

\proclaim{Theorem FB}
For $n\ge1$, the face numbers $f_{k,l}(B_n,m)$ are given by
$$f_{k,l}(B_n,m)={\binom {n} { k+l}} 
{\binom {m n+k-1} k}.$$
Here we identify $B_1$ with $A_1$.
\endproclaim

\demo{Proof}
By inspection, the formula for $f_{k,l}(B_1,m)$ given in the theorem agrees
with the formula for $f_{k,l}(A_1,m)$ in Theorem~FA.
Hence, in view of Proposition~F and (\AD) in Section~2,
in order to prove this claim we have to show
$$l\cdot f_{k,l}(B_n,m)=
\underset l_1+l_2=l-1\to{\underset k_1+k_2=k\hphantom{-1}
\to{\sum _{n_1+n_2=n-1} ^{}}}
f_{k_1,l_1}(B_{n_1},m)\cdot f_{k_2,l_2}(A_{n_2},m)
\tag\BE$$
and
$$
\sum _{k_1+l_1=k} ^{}f_{k_1,l_1}(B_n,m)=
\binom {n}k\binom {mn+k}k.
\tag\BF$$

The triple sum on the right-hand side of (\BE) is 
$$\multline 
\sum _{k_1,l_1,n_1\ge0} ^{}{\binom {n_1} {
k_1+l_1}} 
{\binom {m n_1+k_1-1} {k_1}}\\
\cdot
\frac{l-l_1 }{k-k_1+l-l_1}{\binom
{n-n_1-1} {k- k_1+l-l_1-1}} 
{\binom {m (n-n_1)+k-k_1-1} {k-k_1}}.
\endmultline$$
We replace $n_1$ by $n_1+k_1+l_1$ and rewrite the resulting
expression in the form
$$\multline 
\sum _{k_1,l_1,n_1\ge0} ^{}{\binom {n_1+k_1+l_1} {n_1}} 
{\binom {m (n_1+k_1+l_1)+k_1-1} {k_1}}\\
\cdot
\frac{m(l-l_1) }{m (n-n_1-k_1-l_1)+k-k_1}{\binom
{n-n_1-k_1-l_1} {n-k-l-n_1}}  \\
\cdot
{\binom {m (n-n_1-k_1-l_1)+k-k_1} {k-k_1}}.
\endmultline\tag\BG$$
Forgetting the sum over $l_1$, this is now in the form
of the left-hand side of (\BH) with $n$ replaced by $n-k-l$,
$a=m+1$, $c=m$, $\al=ml_1$, $\al'=m(l-l_1)$, $b=d=1$,
$\be=l_1+1$, $\be'=l-l_1$. Substituting the right-hand side, we obtain
$$
\sum _{l_1=0} ^{l-1}\binom {mn+k-1}{k}
\binom {n}{n-k-l}$$
for the sum in (\BG), or, equivalently,
$$l\cdot{\binom {n} { k+l}} 
{\binom {mn+k-1} k},$$
which is indeed equal to $l\cdot f_{k,l}(B_n,m)$.
This proves (\BE).

\medskip
Next we compute the sum on the left-hand side of (\BF),
$$
\sum _{k_1=0} ^{k}{\binom {n} { k}} 
{\binom {m n+k_1-1} {k_1}}=
{\binom {n} { k}} 
\binom {mn+k} {k},
$$
the simplification of summation being due to the Chu--Vandermonde summation.
This completes the proof.\quad \quad \qed
\enddemo

\subhead 6. The $F$-triangle for $D_n$\endsubhead
The theorem below gives an explicit expression for 
the refined face numbers $f_{k,l}(D_n,m)$, and, thus, of the
$F$-triangle in type $D_n$.

\proclaim{Theorem FD}
For $n\ge2$, the face numbers $f_{k,l}(D_n,m)$ are given by
$$\multline 
f_{k,l}(D_n,m)=
{\binom {n} { k+l}} 
{\binom {m (n-1)+k-1} k}+
m{\binom {n-1} { k+l-1}} 
{\binom {m (n-1)+k-2} {k-1}}\kern1cm\\
-\de_{l,0}\frac {1} {n-1}\binom {n-1}{k-1}\binom {m(n-1)+k-1}k,
\endmultline$$
where $\de_{l,0}$ is the Kronecker delta, that is, it is equal to $1$
if $l=0$, and it is equal to $0$ otherwise.
Here we identify $D_2$ with $A_1^2$, and we identify $D_3$ with $A_3$.
\endproclaim

\vskip.4cm

\demo{Proof}
By inspection, for $n=2$ the formula for $f_{k,l}(D_2,m)$ given in 
the theorem yields for the $F$-triangle
$$F^m_{D_2}(x,y)=
\sum _{k,l\ge0} ^{}f_{k,l}(D_2,m)\,x^ky^l=
m^2 x^2+2 m xy+y^2+2 m x+2 y+1=(m x+y+1)^2,$$
which, according to Theorem~FA, is indeed the $F$-triangle of $A_1^2$.
Furthermore, again by inspection, 
the formula for $f_{k,l}(B_3,m)$ given in the theorem agrees
with the formula for $f_{k,l}(A_3,m)$ in Theorem~FA.
Hence, in view of Proposition~F and (\AD) in Section~2,
in order to prove this claim we have to show
$$l\cdot f_{k,l}(D_n,m)=
\underset l_1+l_2=l-1\to{\underset k_1+k_2=k\hphantom{-1}
\to{\sum _{n_1+n_2=n-1,\ n_1\ge2} ^{}}}
f_{k_1,l_1}(D_{n_1},m)\cdot f_{k_2,l_2}(A_{n_2},m)+
2\cdot f_{k,l-1}(A_{n-1},m)
\tag\BI$$
and
$$
\sum _{k_1+l_1=k} ^{}f_{k_1,l_1}(D_n,m)=
\binom {n}k\binom {m(n-1)+k}k+
\binom {n-2}{k-2}\binom {m(n-1)+k-1}k.
\tag\BJ$$

We start with the proof of (\BI). Clearly, it suffices to consider the case
$l\ge1$ because (\BI) is trivially true for $l=0$. We shall
therefore assume $l\ge1$ from now on.

Using the rewriting
$$\frac{l+1 }{k+l+1}{\binom {n} { k+l}} 
{\binom {m (n+1)+k-1} k}=
\frac{m(l+1) }{m (n+1)+k}{\binom {n+1} { n-k-l}} 
{\binom {m (n+1)+k} k}$$
of the defining expression for $f_{k,l}(A_n,m)$ in Theorem~FA,
the expression on the right-hand side of (\BI) is 
$$\allowdisplaybreaks\align 
\sum _{l_1=0} ^{l-1}
\sum _{k_1=0} ^{k}&
\sum _{n_1=2} ^{n-1}{\binom {n_1} { k_1+l_1}} 
{\binom {m (n_1-1)+k_1-1} {k_1}}\\
&\kern.9cm
\cdot
\frac{m(l-l_1) }{m (n-n_1)+k-k_1}{\binom
{n-n_1} { n-n_1- k+ k_1-l+l_1}} 
{\binom {m (n-n_1)+k-k_1} {k-k_1}}\\\tag\BKa\\
&\kern-.5cm
+m\sum _{l_1=0} ^{l-1}
\sum _{k_1=0} ^{k}
\sum _{n_1=2} ^{n-1}{\binom {n_1-1} { k_1+l_1-1}} 
{\binom {m (n_1-1)+k_1-2} {k_1-1}}\\
&\kern.9cm
\cdot
\frac{m(l-l_1) }{m (n-n_1)+k-k_1}{\binom
{n-n_1} {n-n_1-k+ k_1-l+l_1}} 
{\binom {m (n-n_1)+k-k_1} {k-k_1}}\\\tag\BKb\\
&\kern-.5cm
-
\sum _{k_1=0} ^{k}
\sum _{n_1=2} ^{n-1}
\frac {1} {n_1-1}\binom {n_1-1}{k_1-1}\binom {m(n_1-1)+k_1-1}{k_1}\\
&\kern1.5cm
\cdot
\frac{ml }{m (n-n_1)+k-k_1}{\binom
{n-n_1} {n-n_1-k+ k_1-l}} 
{\binom {m (n-n_1)+k-k_1} {k-k_1}}\\\tag\BKc\\
&\kern-.5cm
+\frac{2l }{k+l}{\binom {n-1} { k+l-1}} 
{\binom {m n+k-1} k}.\tag\BKd
\endalign$$

We treat the three sums in (\BK) separately. We begin with the sum
(\BKa). We extend the sum to all $n_1\ge0$. In order to do so, we must
subtract the terms with $n_1=0$ and $n_1=1$. If $n_1=0$, the summand
is only non-zero for $k_1=l_1=0$ because of the presence of the
binomial $\binom {n_1}{k_1+l_1}$. If $n_1=1$, then,
because of the presence of the binomial
coefficient $\binom {m(n_1-1)+k_1-1}{k_1}=\binom {k_1-1}{k_1}$, 
the summand can be non-zero only if $k_1=0$. 
On the other hand, in that case, the summand can be only
non-zero for $l_1=0$ and $l_1=1$, again because of the presence of the
binomial $\binom {n_1}{k_1+l_1}=\binom 1{l_1}$.
In summary, the sum (\BKa) is equal to
$$\align 
\sum _{l_1=0} ^{l-1}&
\sum _{k_1=0} ^{k}
\sum _{n_1=0} ^{n-1}\Bigg({\binom {n_1} { k_1+l_1}} 
{\binom {m (n_1-1)+k_1-1} {k_1}}\\
&\kern1cm
\cdot
\frac{m(l-l_1) }{m (n-n_1)+k-k_1}{\binom
{n-n_1} {n-n_1-k+ k_1-l+l_1}} 
{\binom {m (n-n_1)+k-k_1} {k-k_1}}\Bigg)\\
&
-\frac{l}{k+l}{\binom
{n-1} {k+l-1}} 
{\binom {m n+k-1} {k}}
-\frac{l}{k+l}{\binom
{n-2} {k+l-1}} 
{\binom {m (n-1)+k-1} {k}}\\
&\kern4cm
-\frac{l-1}{k+l-1}{\binom
{n-2} {k+l-2}} 
{\binom {m (n-1)+k-1} {k}},
\endalign$$
where the second-to-last term corresponds to the summand for
$n_1=k_1=l_1=0$, the next-to-last term corresponds to the summand for
$n_1=1$, $k_1=0$, $l_1=0$, and the last term corresponds to the
summand for $n_1=1$, $k_1=0$, $l_1=1$. In the sum over $n_1,k_1,l_1$,
we replace $n_1$ by $n_1+k_1+l_1$. Forgetting the sum over $l_1$, 
we see that it is then in the form of the left-hand side of (\BH) with 
$n$ replaced by $n-k-l$,
$a=m+1$, $c=m$, $\al=m(l_1-1)$, $\al'=m(l-l_1)$, $b=d=1$,
$\be=l_1+1$, $\be'=l-l_1$. Hence, if we substitute the right-hand side,
the expression simplifies to
$$\multline
\sum _{l_1=0} ^{l-1}{\binom
{n} {k+l}} 
{\binom {m (n-1)+k-1} {k}}
-\frac{l}{k+l}{\binom
{n-1} {k+l-1}} 
{\binom {m n+k-1} {k}}\\
-\frac{l}{k+l}{\binom
{n-2} {k+l-1}} 
{\binom {m (n-1)+k-1} {k}}
-\frac{l-1}{k+l-1}{\binom
{n-2} {k+l-2}} 
{\binom {m (n-1)+k-1} {k}}.\\
\endmultline
\tag\BL
$$
Clearly, the sum over $l_1$ sums the same summand for each $l_1$, so
that the result is that summand multiplied by $l$.

We next turn our attention to the sum (\BKb). The first observation is
that for $k_1=0$ the summand vanishes because of the presence of the
binomial coefficient $\binom {m (n_1-1)+k_1-2} {k_1-1}$. We therefore
replace $k_1$ by $k_1+1$ to obtain
$$\multline 
 m\sum _{l_1=0} ^{l-1}
\sum _{k_1=0} ^{k-1}
\sum _{n_1=2} ^{n-1}{\binom {n_1-1} { k_1+l_1}} 
{\binom {m (n_1-1)+k_1-1} {k_1}}\\
\cdot
\frac{m(l-l_1) }{m (n-n_1)+k-k_1-1}{\binom
{n-n_1} {n-n_1-k+ k_1-l+l_1+1}} 
{\binom {m (n-n_1)+k-k_1-1} {k-k_1-1}}.
\endmultline$$
This time we extend the sum to $n_1\ge1$. In order to do
so, we must subtract the terms with $n_1=1$. In the latter case,
because of the presence of the binomial coefficient
$\binom {n_1-1} { k_1+l_1}$, the
summand will vanish except if $k_1=l_1=0$. Thus, we obtain
$$\multline 
 m\sum _{l_1=0} ^{l-1}
\sum _{k_1=0} ^{k-1}
\sum _{n_1=1} ^{n-1}\Bigg({\binom {n_1-1} { k_1+l_1}} 
{\binom {m (n_1-1)+k_1-1} {k_1}}\\
\cdot
\frac{m(l-l_1) }{m (n-n_1)+k-k_1-1}{\binom
{n-n_1} {n-n_1-k+ k_1-l+l_1+1}} 
{\binom {m (n-n_1)+k-k_1-1} {k-k_1-1}}\Bigg)\\
-\frac{ml }{k+l-1}{\binom
{n-2} {k+l-2}} 
{\binom {m (n-1)+k-2} {k-1}}
\endmultline$$
for (\BKb). In the triple sum,
we replace $n_1$ by $n_1+k_1+l_1+1$. Forgetting the sum over $l_1$, 
we see that it is then in the form of the left-hand side of (\BH) with 
$n$ replaced by $n-k-l$,
$a=m+1$, $c=m$, $\al=ml_1$, $\al'=m(l-l_1)$, $b=d=1$,
$\be=l_1+1$, $\be'=l-l_1$, and $k$ replaced by $k-1$. 
Hence, if we substitute the right-hand side,
the expression simplifies to
$$
m\sum _{l_1=0} ^{l-1}\binom {n-1}{k+l-1}
{\binom {m (n-1)+k-2} {k-1}}
-\frac{ml }{k+l-1}{\binom
{n-2} {k+l-2}} 
{\binom {m (n-1)+k-2} {k-1}}.
\tag\BM$$
Also here, the sum over $l_1$ sums the same summand for each $l_1$, so
that the result is that summand multiplied by $l$.

Finally we treat the sum (\BKc). Again, we want to extend the sum over
$n_1$ to $n_1\ge0$. In order to do so, we would have to 
subtract the terms for $n_1=1$ and $n_1=0$. However, 
it is somewhat unclear which values we should give the summand 
for these choices of $n_1$.
To obtain a partial answer, we rewrite 
$$\frac {1} {n_1-1}\binom {n_1-1}{k_1-1}\binom {m(n_1-1)+k_1-1}{k_1}=
\frac {m} {m(n_1-1)+k_1}\binom {n_1-1}{n_1-k_1}\binom
{m(n_1-1)+k_1}{k_1}.\tag\BMa$$
(This rewriting is already in the spirit of the forth-coming
application of Carlitz's identity (\BD). We alert the reader
that, according to our convention (\BDa), 
the rewriting $\binom {n_1-1}{k_1-1}=\binom
{n_1-1}{n_1-k_1}$ is without problem as long as $n_1\ge1$,
which is the case in (\BKc). However, it
becomes wrong if $1>n_1\ge k_1$, in which case $\binom
{n_1-1}{k_1-1}=0$ while $\binom
{n_1-1}{n_1-k_1}\ne0$. In the following considerations, whenever
we talk about cases where $1>n_1\ge k_1$, we shall talk about
the {\it right-hand side} in (\BMa).)
If $n_1=0$, then, because
of the presence of the binomial coefficient $\binom {n_1-1}{n_1-k_1}$, the
summand is only non-zero if $k_1=0$, in which case it equals 
$$-\frac {ml} {mn+k}\binom n{n-k-l}\binom {mn+k}k=
-\frac {l} {k+l}\binom {n-1}{k+l-1}\binom {mn+k-1}k.$$
For $n_1=1$, the above expression vanishes certainly if $k_1>1$. If
$k_1=1$, it is equal to $m$. But if $k_1=0$, it is still not clear
which value to assign to it. Leaving this question open for the
moment, the arguments so far show that the expression (\BKc) is equal to
$$\multline -
\sum _{k_1=0} ^{k}
\sum _{n_1=0} ^{n-1}\Bigg(
\frac {m} {m(n_1-1)+k_1}\binom {n_1-1}{n_1-k_1}\binom {m(n_1-1)+k_1}{k_1}\\
\cdot
\frac{ml }{m (n-n_1)+k-k_1}{\binom
{n-n_1} {n-n_1-k+ k_1-l}} 
{\binom {m (n-n_1)+k-k_1} {k-k_1}}\Bigg)\\
-\frac{l }{k+l}{\binom
{n-1} {k+l-1}} 
{\binom {m n+k-1} {k}}
+\frac{ml }{k+l-1}{\binom
{n-2} {k+l-2}} 
{\binom {m (n-1)+k-2} {k-1}}\\
+(\text{summand for $n_1=1$, $k_1=0$}).
\endmultline$$
We now replace $n_1$ by $n_1+k_1$ in the double sum. This leads to the
expression
$$\multline -
\sum _{n_1,k_1\ge0} ^{}
\Bigg(
\frac {m} {m(n_1+k_1-1)+k_1}\binom {n_1+k_1-1}{n_1}
\binom {m(n_1+k_1-1)+k_1}{k_1}\\
\cdot
\frac{ml }{m (n-n_1-k_1)+k-k_1}{\binom
{n-n_1-k_1} {n-n_1-k-l}} 
{\binom {m (n-n_1-k_1)+k-k_1} {k-k_1}}\Bigg)\\
-\frac{l }{k+l}{\binom
{n-1} {k+l-1}} 
{\binom {m n+k-1} {k}}
+\frac{ml }{k+l-1}{\binom
{n-2} {k+l-2}} 
{\binom {m (n-1)+k-2} {k-1}}\\
+(\text{summand for $n_1=1$, $k_1=0$}).
\endmultline$$
The double sum is now exactly equal to the negative of 
the left-hand side of (\BD) with $n$ replaced by $n-k-l$,
$a=m+1$, $c=m$, $\al=-m$, $\al'=ml$, $b=d=1$,
$\be=-1$, $\be'=l$. From there, we can also determine the missing
value of the summand for $n_1=1$ and $k_1=0$. Namely, we have
$$A_{0,1}(\al,\be)=\be=-1.$$
Thus, if we substitute the right-hand side of (\BD), we obtain
$$\multline 
\frac {l-1} {k+l-1}\binom {n-2}{k+l-2}\binom {m(n-1)+k-1}k
-\frac{l }{k+l}{\binom
{n-1} {k+l-1}} 
{\binom {m n+k-1} {k}}\\
+\frac{ml }{k+l-1}{\binom
{n-2} {k+l-2}} 
{\binom {m (n-1)+k-2} {k-1}}
+\frac{l }{k+l}{\binom
{n-2} {k+l-1}} 
{\binom {m (n-1)+k-1} {k}}\\
\endmultline\tag\BN$$
for (\BKc). 

Adding the expressions (\BL), (\BM), (\BN) and (\BKd),
we obtain that the sum in (\BK) is equal to
$$l{\binom {n} { k+l}} 
{\binom {m(n-1)+k-1} k}+
ml{\binom {n-1} { k+l-1}} 
{\binom {m(n-1)+k-2} {k-1}},$$
which is indeed equal to $l\cdot f_{k,l}(D_n,m)$ if $l\ge1$.
This proves (\BI).

\medskip
Next we compute the sum on the left-hand side of (\BJ),
$$\align \sum _{k_1=0} ^{k}{\binom {n} { k}} 
&{\binom {m (n-1)+k_1-1} {k_1}}+
m\sum _{k_1=0} ^{k}{\binom {n-1} { k-1}} 
{\binom {m (n-1)+k_1-2} {k_1-1}}\\
&\kern2cm
-\frac {1} {n-1}\binom {n-1}{k-1}\binom {m(n-1)+k-1}k\\
&={\binom {n} { k}} 
{\binom {m (n-1)+k} {k}}+
m{\binom {n-1} { k-1}} 
{\binom {m (n-1)+k-1} {k-1}}\\
&\kern2cm
-\frac {1} {n-1}\binom {n-1}{k-1}\binom {m(n-1)+k-1}k\\
&={\binom {n} { k}} 
{\binom {m (n-1)+k} {k}}+
{\binom {n-2} { k-2}} 
{\binom {m (n-1)+k-1} {k}},
\endalign$$
the simplification of summations being due to the Chu--Vandermonde summation.
This completes the proof.\quad \quad \qed
\enddemo

\subhead 7. The $F$-triangle in the exceptional cases\endsubhead
It is a routine matter to use Proposition~F in Section~2 (and a
computer algebra package) to find the
$F$-triangles for the exceptional root systems. We list our findings
below. 

\medskip\noindent
{\it The $F$-triangle for $I_2(a)$}: 
$$F^m_{I_2(a)}(x,y)= \frac{m (m a+a-2)}{2}  x^2+2 m  xy+a m x+y^2+2 y+1.
\tag\CA$$

\medskip\noindent
{\it The $F$-triangle for $H_3$}: 
$$\multline 
F^m_{H_3}(x,y)= \frac{m (5 m+2) (5 m+4)}{3}  x^3+m (5 m+2)  x^2y\\
+
5 m (5 m+2) x^2+3 m xy^2
    +10 m  xy+15 m x+y^3+3 y^2+3 y+1.
\endmultline\tag\CB$$

\medskip\noindent
{\it The $F$-triangle for $H_4$}: 
$$\multline 
F^m_{H_4}(x,y)=  \frac{ m (3 m+1) (5 m+3) (15 m+14)}{4} x^4+m (3 
    m+1) (5 m+3) x^3y\\
   +15 m (3
m+1) (5 m+3) x^3+\frac{1}{2} m (17 m+5)  x^2y^2+m (45 m+14) x^2y\\
+\frac{1}{2} m (465
    m+149) x^2
   +4 m xy^3+17 m xy^2+31 m xy+60 m x+y^4+4
    y^3+6 y^2+4 y+1.
\endmultline\tag\CC$$

\medskip\noindent
{\it The $F$-triangle for $F_4$}: 
$$\multline 
F^m_{F_4}(x,y)=  \frac{m (2 m+1) (3 m+1) (6 m+5)}{2}  x^4
+2 m (2
    m+1) (3 m+1) x^3y\\+12 m (2 m+1) (3 m+1) x^3+2 m (4 m+1) x^2 y^2 +2 m (18 m+5) 
    x^2y+m (78 m+23) x^2\\
+4 m xy^3+16 m xy^2+26 m xy+24 m x+y^4+4 y^3+6 y^2+4 y+1.
\endmultline\tag\CD$$

\medskip\noindent
{\it The $F$-triangle for $E_6$}: 
$$\allowdisplaybreaks
\multline 
F^m_{E_6}(x,y)=  \frac{1}{30} m (2 m+1) (3 m+1) (4 m+1) (6 m+5) (12 m+7) x^6\\
+\frac{1}{5} m (2 m+1) (3 m+1) (4 m+1)
    (12 m+7)  x^5y+\frac{6}{5} m
    (2 m+1) (3 m+1) (4 m+1) (12 m+7) x^5\\
+\frac{1}{2} m (3 m+1) (4 m+1) (8 m+3)  x^4y^2
+2 m (3 m+1) (4 m+1) (12 m+5)  x^4y\\
+2 m (3 m+1)
    (4 m+1) (30 m+13) x^4+\frac{5}{3} m (4
    m+1) (5 m+1)  x^3y^3\\
+m (4 m+1) (48 m+11)  x^3y^2+m (4 m+1) (120 m+31)  x^3y
+9 m (4 m+1) (18 m+5)
    x^3\\
+\frac{5}{2} m (7 m+1)  x^2y^4+5 m (20 m+3)
     x^2y^3+m (242 m+39)  x^2y^2+3 m (108 m+19)  x^2y\\
+12 m (21 m+4) x^2+6 m
     xy^5+35 m  xy^4+85 m  xy^3+111 m  xy^2+84 m xy+36 m x\\
+y^6+6 y^5+15
    y^4+20 y^3+15 y^2+6 y+1. 
\endmultline\tag\CE$$

\medskip\noindent
{\it The $F$-triangle for $E_7$}: 
{\eightpoint
$$\allowdisplaybreaks
\multline 
F^m_{E_7}(x,y)=  \frac{1}{280} m (3 m+1) (3 m+2) (9 m+2) (9 m+4) (9 m+5) (9 m+8)
    x^7\\
+\frac{9}{40} m (3 m+1) (3 m+2) (9 m+2) (9 m+4) (9 m+5)
    x^6
+\frac{1}{40} m (3 m+1) (3 m+2) (9 m+2) (9 m+4) (9 m+5) 
    x^6y\\
+\frac{3}{40} m (3 m+1) (7 m+3) (9 m+2) (9 m+4)  x^5y^2
+\frac{3}{20} m (3 m+1) (9 m+2) (9
    m+4) (27 m+13)  x^5y\\
+\frac{3}{40} m
    (3 m+1) (9 m+2) (9 m+4) (207 m+103) x^5
+\frac{1}{8} m (3 m+1) (9 m+2) (27 m+7) 
    x^4y^3\\
+\frac{3}{8} m (3 m+1) (9 m+2) (63 m+19)  x^4y^2
+\frac{3}{8} m (3 m+1) (9 m+2) (207 m+71)
 x^4y\\
+\frac{21}{8} m (3 m+1)
    (9 m+2) (63 m+23) x^4+m (6
    m+1) (9 m+2)  x^3y^4+\frac{3}{2} m (9 m+2) (27 m+5)  x^3y^3\\
+\frac{3}{2} m
    (9 m+2) (81 m+17)  x^3y^2+\frac{21}{2}
    m (9 m+2) (21 m+5)  x^3y
+\frac{21}{2} m (9 m+2) (27 m+7) x^3
+3 m (8 m+1)  x^2y^5\\
+3 m (54 m+7) 
    x^2y^4+\frac{3}{2} m (315 m+43)  x^2y^3+\frac{21}{2} m (75 m+11) 
    x^2y^2+\frac{21}{2} m (81 m+13)  x^2y
  +\frac{21}{2} m (63 m+11) x^2\\
+7 m 
    xy^6+48 m  xy^5+141 m  xy^4+231 m  xy^3+231 m  xy^2+147 m  xy+63 m x\\+y^7+7
    y^6+21 y^5+35 y^4+35 y^3+21 y^2+7 y+1.
\endmultline\tag\tenpoint\CF$$}

\medskip\noindent
{\it The $F$-triangle for $E_8$}: 
{\eightpoint
$$\allowdisplaybreaks
\multline 
F^m_{E_8}(x,y)
=   \frac{m (3 m+1) (5 m+1) (5 m+2) (5 m+3) (15 m+8) (15 m+11) (15 m+14)
    x^8}{1344}\\
+\frac{1}{168} m (3 m+1) (5 m+1) (5 m+2) (5 m+3) (15 m+8) (15
    m+11) x^7y\\
+\frac{5}{56} m (3 m+1) (5 m+1) (5 m+2) (5 m+3) (15 m+8) (15
    m+11) x^7\\
+\frac{1}{48} m (3 m+1) (5 m+1) (5 m+2) (15 m+7) (15 m+8) 
    x^6y^2\\
+\frac{5}{24} m (3 m+1) (5 m+1) (5 m+2) (15 m+8)^2 x^6y
+\frac{5}{48} m (3 m+1) (5 m+1) (5 m+2) (15 m+8) (195 m+107)
    x^6\\
+\frac{1}{3} m
    (3 m+1) (5 m+1) (5 m+2) (10 m+3) x^5y^3+\frac{5}{8} m (3 m+1) (5 m+1) (5
    m+2) (45 m+16) x^5y^2\\
+\frac{25}{8} m (3 m+1) (5 m+1) (5 m+2) (39 m+16) x^5y
+15 m (3 m+1) (5 m+1) (5 m+2) (30 m+13) x^5\\
+\frac{1}{6} m
    (5 m+1) (10 m+3) (19 m+4) x^4y^4+m (5 m+1) (10 m+3) (25 m+6) 
    x^4y^3\\
+\frac{1}{4} m (5 m+1) \left(3675 m^2+2125 m+308\right) 
    x^4y^2+m (5 m+1)
    \left(2250 m^2+1395 m+218\right) x^4y\\
+\frac{1}{2} m (5 m+1) \left(10350 m^2+6675 m+1084\right) x^4
+\frac{7}{3} m (5 m+1) (7 m+1) 
    x^3y^5\\
+\frac{1}{3} m (5 m+1) (380 m+59) x^3y^4+\frac{1}{3} m (5 m+1) (1315
    m+226) x^3y^3
+m (5 m+1) (915 m+178) x^3y^2\\
+m (5 m+1) (1380 m+307) x^3y
+45 m (5 m+1) (45 m+11) x^3+\frac{7}{2} m (9 m+1) x^2y^6+7 m (35 m+4) 
    x^2y^5\\
+\frac{1}{2} m (1675 m+199) x^2y^4+4 m (415 m+52) x^2y^3+\frac{1}{2}
    m (4295 m+579) x^2y^2+75 m (27 m+4) x^2y\\
+\frac{35}{2} m (105 m+17) x^2+8
    m xy^7+63 m xy^6+217 m xy^5+428 m xy^4+532 m xy^3+435 m xy^2+245 m xy+120 m
    x\\
+y^8+8 y^7+28 y^6+56 y^5+70 y^4+56 y^3+28 y^2+8 y+1.
\endmultline\tag\tenpoint\CG$$}

\subhead 8. The $F=M$ Conjecture\endsubhead
In order to state the ``$F=M$ Conjecture" for generalised cluster
complexes, we need to first introduce Armstrong's \cite{\ArmDAA}
$m$-divisible non-crossing partition posets. 

Given a root system $\Phi$ and an element $\al\in\Phi$, let $t_\al$
denote the corresponding reflection in the central hyperplane perpendicular to
$\al$. Let $W=W(\Phi)$ be the group generated by these reflections. By
definition, any element $w$ of $W$ can be represented as a product 
$w=t_1t_2\cdots t_\ell$, where the $t_i$'s are reflections. We call
the minimal number of reflections which is needed for such a
product representation the {\it absolute length\/} of $w$, and we
denote it by $\ell_T(w)$. We then define the {\it absolute order} on
$W$, denoted by $\le_T$, by 
$$u\le_T w\quad \text{if and only if}\quad
\ell_T(w)=\ell_T(u)+\ell_T(u^{-1}w).$$
It can be shown that this is equivalent to the statement that any
shortest product representation of $u$ by reflections
occurs as an initial segment in some shortest product representation
of $w$ by reflections. 

We can now define the {\it non-crossing partition lattice
$NC(\Phi)$}. Let $c$ be a {\it Coxeter element\/} in $W$, that is, the
product of all reflections corresponding to the simple roots. Then
$NC(\Phi)$ is defined to be
the restriction of the partial order $\le_T$ to the set
of all elements which are less than or equal to $c$ in absolute order. 
This definition makes sense because, regardless
of the chosen Coxeter element $c$, the resulting poset
is always the same up to isomorphism.
It can be shown that $NC(\Phi)$ is indeed a lattice. (See
\cite{\BRWaAB} for a uniform proof.)
The term ``non-crossing partition lattice" is used because
$NC(A_n)$ is isomorphic to the lattice of non-crossing partitions
originally introduced by Kreweras \cite{\KrewAC}, and because also
$NC(B_n)$ and $NC(D_n)$ can be realized as lattices of
non-crossing partitions (see \cite{\AtReAA, \ReivAG}).

The poset of {\it $m$-divisible non-crossing partitions} has as a
groundset the following subset of $(NC(\Phi))^{m+1}$,
$$\multline
NC^m(\Phi)=\{(w_0;w_1,\dots,w_m):w_0w_1\cdots w_m=c\text{ and }\\
\ell_T(w_0)+\ell_T(w_1)+\dots+\ell_T(w_m)=\ell_T(c)\}.
\endmultline\tag\CH$$
The order relation is defined by
$$(u_0;u_1,\dots,u_m)\le(w_0;w_1,\dots,w_m)\quad \text{if and only
if}\quad u_i\ge_T w_i,\ 1\le i\le m.$$
We emphasize that, according to this definition, $u_0$ and $w_0$ need
not be related in any way. The poset $NC^m(\Phi)$ is graded by 
the rank function
$$\rk\big((w_0;w_1,\dots,w_m)\big)=\ell_T(w_0).$$
Thus, there is a unique maximal element, namely $(c;\ep,\dots,\ep)$,
where $\ep$ stands for the identity element in $W$, but, if $m>1$, there 
are many
different minimal elements. In particular, there is no global minimum
in $NC^m(\Phi)$ if $m>1$ and, hence, $NC^m(\Phi)$ is not a lattice
for $m>1$. (It is, however, a graded join-semilattice, see
\cite{\ArmDAA, Theorem~2.2.7}.) 
We remark that for $NC^m(A_n)$ and $NC^m(B_n)$
combinatorial realisations are available as subposets of
non-crossing partitions in which each block has a size which is
divisible by $m$. The corresponding translations are due to Armstrong
\cite{\ArmDAA, Sec.~3}. In type $A_n$, the resulting poset had been
earlier studied by Edelman \cite{\EdelAA}. The analogous combinatorial
realisation of $NC^m(D_n)$, generalising the one of Athanasiadis and
Reiner \cite{\AtReAA} for $m=1$, has not yet been worked out.

Next, we define the ``$M$-triangle" of $NC^m(\Phi)$ as
$$
M^m_\Phi(x,y)=\sum _{u,w\in NC^m(\Phi)} ^{}\mu(u,w)\,x^{\rk u}y^{\rk w},$$
where $\mu(u,w)$ is the M\"obius function in $NC^m(\Phi)$.

The generalised version of Chapoton's (ex-)conjecture \cite{\ChaFAA,
Conjecture~1}, due to Armstrong \cite{\ArmDAA, Sec.~4}, is the following.

\proclaim{Conjecture FM}
For any finite root system $\Phi$ of rank $n$, we have
$$F^m_\Phi(x,y)=y^n\,M^m_\Phi\(\frac {1+y} {y-x},\frac {y-x} {y}\).$$
Equivalently, 
$$(1-xy)^n F^m_\Phi\left(\frac {x(1+y)} {1-xy},\frac {xy}
{1-xy}\right)=
\sum _{u,w\in(NC^m(\Phi))^*}
^{}\mu^*(u,w)\,(-x)^{\rk^*w}(-y)^{\rk^*u},\tag\BO$$ 
where $(NC^m(\Phi))^*$ denotes the poset {\rm dual} to $NC^m(\Phi)$
{\rm(}i.e., the poset which arises from $NC^m(\Phi)$ by reversing all
order relations{\rm)}, where $\mu^*$ denotes the M\"obius
function in\linebreak 
$(NC^m(\Phi))^*$, and where $\rk^*$ denotes the rank
function in $(NC^m(\Phi))^*$.
\endproclaim

Since the M\"obius function is multiplicative (see e.g.\
\cite{\StanAP, Prop.~3.8.2}), the multiplicativity property (\AA)
for the $F$-triangle holds also for the $M$-triangle.
Therefore, it is enough to prove the conjecture for the irreducible
root systems. In Sections~9--17, we provide proofs for the root
systems of type $A_n$, $B_n$, $I_2(a)$, $H_3$, $H_4$, $F_4$,
and $E_6$,
and a partial proof for the root system of type $D_n$.

\subhead 9. Proof of the $F=M$ Conjecture for $A_n$\endsubhead
In this section we prove Conjecture~FM for the type $A_n$. 
In the spirit of this paper, we follow a computational approach. 
We first simplify the left-hand side of (\BO) by a double application
of the Chu--Vandermonde summation. Subsequently, we compute the
right-hand side of (\BO) by relying on a result on
rank selected chain enumeration in the $m$-divisible non-crossing
partition lattice in type $A_n$ due to Edelman \cite{\EdelAA}. 

The link between chain enumeration and the M\"obius function is the
following. (The reader should consult \cite{\StanAP, Sec.~3.11} for
more information on this topic.) 
Given a poset $P$ and two elements $u$ and $w$, $u\le w$,
in the poset, the {\it zeta polynomial\/} of the interval $[u,w]$,
denoted by $Z(u,w;z)$, is the number of (multi)chains from $u$ to
$w$ of length $z$. (It can be shown that this is indeed a
polynomial in $z$.) Then the M\"obius function of $u$ and $w$ is equal
to $\mu(u,w)=Z(u,w;-1)$.

\proclaim{Proposition A}
In type $A_n$, the left-hand side of {\rm(\BO)} is equal to
$$
\sum _{r,s\ge0} ^{}x^sy^r\frac {1} {s+1}\binom ns\binom {m(n+1)}r
\binom {m(n+1)+s-r-1}{s-r}.\tag\BP$$
\endproclaim
\demo{Proof}
By definition of $F^m_{A_n}(x,y)$, and by Theorem~FA
in Section~4, the left-hand side of
(\BO) in type~$A_n$ is equal to
$$
\sum _{k,l,r,s\ge0} ^{}\frac{l+1 }{k+l+1}{\binom {n} { k+l}}
{\binom {m (n+1)+k-1} k}\binom kr\binom {n-k-l}s(-1)^s
x^{k+l+s}y^{l+r+s}.$$
Fixing $L=s+l$, we rewrite this as
$$\multline
\sum _{k,L,r\ge0} ^{}
\frac {n!\,(m (n+1)+k-1)!} {(m (n+1)-1)!\,r!\,(k-r)!\,(n-k-L)!\,(k+L+1)!}
x^{k+L}y^{r+L}\\
\cdot\sum _{s=0} ^{L}{(L-s+1)} (-1)^s\binom{k+L+1}s.
\endmultline$$
We compute the sum over $s$ by the Chu--Vandermonde
summation. Thus, we arrive at
$$
\sum _{k,L,r\ge0} ^{}
\frac {n!\,(m (n+1)+k-1)!} {(m (n+1)-1)!\,r!\,(k-r)!\,(n-k-L)!\,(k+L+1)!}
x^{k+L}y^{r+L}
(-1)^L\binom{k+L-1}L.
$$
We now write $K=k+L$ and $R=r+L$.
Subsequently, the sum over $L$ can be computed using the Chu--Vandermonde
summation. The result is
$$
\sum _{K,R\ge0} ^{}
\frac {n!\,(m (n+1)+K-R-1)!\,(m(n+1))!} 
{(m (n+1)-1)!\,(m (n+1)-R)!\,R!\,(K-R)!\,(n-K)!\,(K+1)!}
x^{K}y^{R}.
$$
Aside from a parameter replacement, this is exactly the expression
(\BP).\quad \quad \qed
\enddemo

For the computation of the right-hand side of (\BO) we require 
the following theorem due to Edelman
\cite{\EdelAA}.


\proclaim{Theorem NA}
The number of chains in $(NC^m(A_n))^*$ with successive rank jumps
$s_1,s_2,\dots,s_\ell$, $s_1+s_2+\dots+s_\ell=n$, is
$$\binom {m(n+1)}{s_1}\cdots \binom {m(n+1)}{s_{\ell-1}}
\frac {1} {n+1}\binom {n+1}{s_\ell}.\tag\BPa$$
\endproclaim

\demo{Proof of Conjecture~FM in type $A_n$}
We now compute the right-hand side of (\BO), that is,
$$
\sum _{u,w\in 
(NC^m(A_n))^*} ^{}\mu^*(u,w)(-x)^{\rk^* w}(-y)^{\rk^* u}.$$
In order to compute the coefficient of $x^sy^r$ in this expression,
$$
(-1)^{r+s}\underset\text{ with }\rk^* u=r\text{ and }\rk^* w=s\to
{\sum _{u,w\in (NC^m(A_n))^*} ^{}}\mu^*(u,w),$$
we compute the sum of all corresponding zeta polynomials (in the
variable $z$), multiplied by $(-1)^{r+s}$,
$$
(-1)^{r+s}\underset\text{ with }\rk^* u=r\text{ and }\rk^* w=s\to
{\sum _{u,w\in (NC^m(A_n))^*} ^{}}Z(u,w;z),
$$
and then put $z=-1$. 

For computing this sum of zeta polynomials,
we must set
$\ell=z+2$, $s_1=r$, $n-s_\ell=s$, $s_2+s_3+\dots+s_{\ell-1}=s-r$ in (\BPa), and
then sum the resulting expression over all possible
$s_2,s_3,\dots,s_{\ell-1}$. By using the Chu--Vandermonde summation,
one obtains
$$\frac {1} {n+1}\binom {m(n+1)}{r}\binom {zm(n+1)}{s-r}
\binom {n+1}{n-s}.$$
If we put $z=-1$ in this expression and multiply it by $(-1)^{r+s}$,
then we obtain exactly the coefficient of $x^sy^r$ in (\BP).\quad \quad
\qed
\enddemo

\subhead 10. Proof of the $F=M$ Conjecture for $B_n$\endsubhead
In this section we prove Conjecture~FM for the type $B_n$,
by following the same approach as the one for type $A_n$ in the
previous section.

\proclaim{Proposition B}
In type $B_n$, the left-hand side of {\rm(\BO)} is equal to
$$
\sum _{r,s\ge0} ^{}x^sy^r\binom ns\binom {mn}r
\binom {mn+s-r-1}{s-r}.\tag\BQ$$
\endproclaim
\demo{Proof}
By definition of $F^m_{B_n}(x,y)$, and by Theorem~FB
in Section~5, the left-hand side of
(\BO) in type~$B_n$ is equal to
$$
\sum _{k,l,r,s\ge0} ^{}{\binom {n} { k+l}}
{\binom {m n+k-1} k}\binom kr\binom {n-k-l}s(-1)^s
x^{k+l+s}y^{l+r+s}.$$
Fixing $L=s+l$, we rewrite this as
$$
\sum _{k,L,r\ge0} ^{}
\frac {n!\,(m n+k-1)!} {(m n-1)!\,r!\,(k-r)!\,(n-k-L)!\,(k+L)!}
x^{k+L}y^{L+r}\sum _{s=0} ^{L} (-1)^s\binom{k+L}s.
$$
We compute the sum over $s$ by the Chu--Vandermonde
summation. Thus, we arrive at
$$
\sum _{k,L,r\ge0} ^{}
\frac {n!\,(m n+k-1)!} {(m n-1)!\,r!\,(k-r)!\,(n-k-L)!\,(k+L)!}
x^{k+L}y^{L+r}
(-1)^L\binom{k+L-1}L.
$$
We now write $K=k+L$ and $R=r+L$.
Subsequently, the sum over $L$ can be computed using the Chu--Vandermonde
summation. The result is
$$
\sum _{K,R\ge0} ^{}
\frac {n!\,(m n+K-R-1)!\,(mn)!} 
{(m n-1)!\,(m n-R)!\,R!\,(K-R)!\,(n-K)!\,K!}
x^{K}y^{R}.
$$
Aside from a parameter replacement, this is exactly the expression
(\BQ).\quad \quad \qed
\enddemo

For the computation of the right-hand side of (\BO) we require 
the following theorem due to Armstrong
\cite{\ArmDAA, Theorem~3.5.7}.


\proclaim{Theorem NB}
The number of chains in $(NC^m(B_n))^*$ with successive rank jumps
$s_1,s_2,\dots,s_\ell$, $s_1+s_2+\dots+s_\ell=n$, is
$$\binom {mn}{s_1}\cdots \binom {mn}{s_{\ell-1}}
\binom {n}{s_\ell}.\tag\BQa$$
\endproclaim

\demo{Proof of Conjecture~FM in type $B_n$}
We now compute the right-hand side of (\BO), that is,
$$
\sum _{u,w\in (NC^m(B_n))^*} ^{}\mu^*(u,w)
(-x)^{\rk^* w}(-y)^{\rk^* u}.$$
In order to compute the coefficient of $x^sy^r$ in this expression,
$$
(-1)^{r+s}\underset\text{ with }\rk^* u=r\text{ and }\rk^* w=s\to
{\sum _{u,w\in (NC^m(B_n))^*} ^{}}\mu^*(u,w),$$
we compute the sum of all corresponding zeta polynomials (in the
variable $z$), multiplied by $(-1)^{r+s}$,
$$
(-1)^{r+s}\underset\text{ with }\rk^* u=r\text{ and }\rk^* w=s\to
{\sum _{u,w\in (NC^m(B_n))^*} ^{}}Z(u,w;z),
$$
and then put $z=-1$. 

For computing this sum of zeta polynomials,
we must set
$\ell=z+2$, $s_1=r$, $n-s_\ell=s$, $s_2+s_3+\dots+s_{\ell-1}=s-r$ in (\BQa), and
then sum the resulting expression over all possible
$s_2,s_3,\dots,s_{\ell-1}$. By using the Chu--Vandermonde summation,
one obtains
$$\binom {mn}{r}\binom {zmn}{s-r}
\binom {n}{n-s}.$$
If we put $z=-1$ in this expression and multiply it by $(-1)^{r+s}$,
then we obtain exactly the coefficient of $x^sy^r$ in (\BQ).\quad \quad
\qed
\enddemo

\subhead 11. Towards a proof of the $F=M$ Conjecture for $D_n$\endsubhead
This section exhibits how far the approach of the previous two sections
of proving Conjecture~FM in types $A_n$ and $B_n$ can take us in type
$D_n$. The simplification of the left-hand side of (\BO) along the
lines of the proofs of Propositions~A and B goes through smoothly. The
problem which we face in type $D_n$ is that, up to this date, the rank
selected chain enumeration result for $NC^m(D_n)$ has not been found
yet. Thus,
we do not have the means to compute the $M$-triangle for $NC^m(D_n)$. The only
exception is for $m=1$. Namely, for the (ordinary) non-crossing partition
lattice $NC(D_n)=NC^1(D_n)$, Athanasiadis and Reiner \cite{\AtReAA} have
done the rank selected chain enumeration as we need it in our
application. Hence, we are able to prove Conjecture~FM in type $D_n$
if $m=1$.

\proclaim{Proposition D}
In type $D_n$, the left-hand side of {\rm(\BO)} is equal to
$$\multline
\sum _{r,s\ge0} ^{}x^sy^r\Bigg(
2\binom {n-1}{s-1}\binom {m(n-1)}r\binom {m(n-1)+s-r-1}{s-r}\\+
\binom {n-2}{s}\binom {m(n-1)}r\binom
{m(n-1)+s-r-1}{s-r}\\
+m\binom {n-1}{s-1}
\binom {m(n-1)-1}{r-2}\binom {m(n-1)+s-r-1}{s-r}\\-
m\binom {n-1}{s-1}\binom {m(n-1)}r\binom {m(n-1)+s-r-2}{s-r-2}
\Bigg).
\endmultline\tag\BR$$
\endproclaim
\demo{Proof}
By definition of $F^m_{D_n}(x,y)$, and by Theorem~FD
in Section~6, the left-hand side of
(\BO) in type~$D_n$ is equal to
$$\align
\sum _{k,l,r,s\ge0} ^{}&{\binom {n} { k+l}}
{\binom {m (n-1)+k-1} k}\binom kr\binom {n-k-l}s(-1)^s
x^{k+l+s}y^{l+r+s}\tag\BSa\\
&+m\sum _{k,l,r,s\ge0} ^{}{\binom {n-1} { k+l-1}}
{\binom {m (n-1)+k-2} {k-1}}\binom kr\binom {n-k-l}s(-1)^s
x^{k+l+s}y^{l+r+s}\tag\BSb\\
&-\frac {1} {n-1}\sum _{k,r,s\ge0} ^{}{\binom {n-1} { k-1}}
{\binom {m (n-1)+k-1} k}\binom kr\binom {n-k}s(-1)^s
x^{k+s}y^{r+s}.
\tag\BSc\endalign$$
We treat the three sums in (\BS) separately. We begin with the sum
(\BSa). Fixing $L=s+l$, we rewrite it as
$$
\sum _{k,L,r\ge0} ^{}
\frac {n!\,(m (n-1)+k-1)!} {(m (n-1)-1)!\,r!\,(k-r)!\,(n-k-L)!\,(k+L)!}
x^{k+L}y^{L+r}\sum _{s=0} ^{L} (-1)^s\binom{k+L}s.
$$
We compute the sum over $s$ by the Chu--Vandermonde
summation. Thus, we arrive at
$$
\sum _{k,L,r\ge0} ^{}
\frac {n!\,(m (n-1)+k-1)!} {(m (n-1)-1)!\,r!\,(k-r)!\,(n-k-L)!\,(k+L)!}
x^{k+L}y^{L+r}
(-1)^L\binom{k+L-1}L.
$$
We now write $K=k+L$ and $R=r+L$.
Subsequently, the sum over $L$ can be computed using the Chu--Vandermonde
summation. The result is
$$
\sum _{K,R\ge0} ^{}
\frac {n!\,(m (n-1)+K-R-1)!\,(m(n-1))!} 
{(m (n-1)-1)!\,(m (n-1)-R)!\,R!\,(K-R)!\,(n-K)!\,K!}
x^{K}y^{R}.
\tag\BT$$

Next we consider the sum (\BSb). Again fixing $L=s+l$, we rewrite it as
$$\multline
\sum _{k,L,r\ge0} ^{}
\frac {mk\,(n-1)!\,(m (n-1)+k-2)!} {(m
(n-1)-1)!\,r!\,(k-r)!\,(n-k-L)!\,(k+L-1)!}
x^{k+L}y^{L+r}\\
\cdot\sum _{s=0} ^{L} (-1)^s\binom{k+L-1}s.
\endmultline$$
We compute the sum over $s$ by the Chu--Vandermonde
summation. Thus, we arrive at
$$
\sum _{k,L,r\ge0} ^{}
\frac {mk\,(n-1)!\,(m (n-1)+k-2)!} {(m
(n-1)-1)!\,r!\,(k-r)!\,(n-k-L)!\,(k+L-1)!}
x^{k+L}y^{L+r}
(-1)^L\binom{k+L-2}L.
$$
We now write $K=k+L$ and $R=r+L$. Because of the presence of the
factor $k$ in the numerator, this makes a factor of $K-L$ appear.
We split the sum into two parts accordingly, and then, in both parts,
the sum over $L$ can be computed using the Chu--Vandermonde
summation. The result is
$$\multline
\sum _{K,R\ge0} ^{}
\frac {mK\,(n-1)!\,(m (n-1)+K-R-2)!\,(m(n-1))!} 
{(m (n-1)-1)!\,(m (n-1)-R)!\,R!\,(K-R)!\,(n-K)!\,(K-1)!}
x^{K}y^{R}\\-
\sum _{K,R\ge0} ^{}
\frac {m(K-2)\,(n-1)!\,(m (n-1)+K-R-2)!\,(m(n-1))!} 
{(m (n-1)-1)!\,(m (n-1)-R+1)!\,(R-1)!\,(K-R)!\,(n-K)!\,(K-1)!}
x^{K}y^{R}.
\endmultline\tag\BU$$

Finally, we turn to the sum (\BSc).
We write $K=k+s$ and $R=r+s$.
Subsequently, the sum over $s$ can be computed using the Chu--Vandermonde
summation. The result is
$$
-\sum _{K,R\ge0} ^{}\frac {1} {n-1}
\frac {(n-1)!\,(m (n-1)+K-R-1)!\,(m(n-1))!} 
{(m (n-1)-1)!\,(m (n-1)-R)!\,R!\,(K-R)!\,(n-K)!\,(K-1)!}
x^{K}y^{R}.\tag\BV
$$

To complete the proof, we add the expressions (\BT), (\BU) and (\BV).
Doing the parameter replacements $K\to s$, $R\to r$ and minor rewriting, 
this leads to the expression
$$\multline
\sum _{r,s\ge0} ^{}x^sy^r\Bigg(
\binom ns\binom {m(n-1)}r\binom {m(n-1)+s-r-1}{s-r}\\+
\frac {ms} {s-r}\binom {n-1}{s-1}\binom {m(n-1)}r\binom
{m(n-1)+s-r-2}{s-r-1}\\
+
\frac {m(s-2)} {s-r}\binom {n-1}{s-1}
\binom {m(n-1)}{r-1}\binom {m(n-1)+s-r-2}{s-r-1}\\-
\frac {1} {n-1}\binom {n-1}{s-1}\binom {m(n-1)}r\binom {m(n-1)+s-r-1}{s-r}
\Bigg).
\endmultline$$
It is now a routine verification to show that an alternative way to write this 
is (\BR).\quad \quad \qed
\enddemo

For the computation of the right-hand side of (\BO) in the case that
$m=1$, we require 
the following result due to
Athanasiadis and Reiner \cite{\AtReAA, Theorem~1.2(ii)}.

\proclaim{Theorem ND}
The number of chains 
in $NC(D_n)$ with successive rank jumps $s_1,s_2,\mathbreak\dots,s_\ell$,
$s_1+s_2+\dots+s_\ell=n$, is given by
$$2\binom {n-1}{s_1}\cdots\binom {n-1}{s_\ell}+
\sum _{i=1} ^{\ell}\binom {n-1}{s_1}\cdots\binom {n-2}{s_i-2}\cdots
\binom {n-1}{s_\ell}.\tag\BX$$
\endproclaim

\demo{Proof of Conjecture~FM for $m=1$ in type $D_n$}
If $m=1$, then $NC^m(D_n)$ reduces to the ordinary non-crossing
partition lattice $NC(D_n)$, which is self-dual, that is, 
$(NC(D_n))^*=NC(D_n)$.
Hence, the right-hand side of (\BO) with $m=1$ in type $D_n$ is equal to
$$
\sum _{u,w\in NC(D_n)} ^{}\mu(u,w)(-x)^{\rk w}(-y)^{\rk u}.$$
In order to compute the coefficient of $x^sy^r$ in this expression,
$$
(-1)^{r+s}\underset\text{ with }\rk u=r\text{ and }\rk w=s\to
{\sum _{u,w\in NC(D_n)} ^{}}\mu(u,w),$$
we compute the sum of all corresponding zeta polynomials (in the
variable $z$), multiplied by $(-1)^{r+s}$,
$$
(-1)^{r+s}\underset\text{ with }\rk u=r\text{ and }\rk w=s\to
{\sum _{u,w\in NC(D_n)} ^{}}Z(u,w;z),
$$
and then put $z=-1$. 

For computing this sum of zeta polynomials,
we must set
$\ell=z+2$, $s_1=r$, $n-s_\ell=s$, $s_2+s_3+\dots+s_{\ell-1}=s-r$ in (\BX), and
then sum the resulting expression over all possible
$s_2,s_3,\dots,s_{\ell-1}$. By using the Chu--Vandermonde summation again,
one obtains
$$\multline 
2\binom {n-1}{r}\binom{z(n-1)}{s-r}\binom {n-1}{n-s}+
\binom {n-2}{r-2}\binom{z(n-1)}{s-r}\binom {n-1}{n-s}\\+
\binom {n-1}{r}\binom{z(n-1)}{s-r}\binom {n-2}{n-s-2}+
z\binom {n-1}{r}\binom{z(n-1)-1}{s-r-2}\binom {n-1}{n-s}.
\endmultline$$
If we put $z=-1$ in this expression and multiply it by $(-1)^{r+s}$,
then we obtain exactly the coefficient of $x^sy^r$ in (\BR) with
$m=1$.\quad \quad
\qed
\enddemo

\subhead 12. How to prove the $F=M$ Conjecture in the exceptional
cases\endsubhead
While, at first sight, for a given exceptional root system $\Phi$, it
seems that computing the $M$-triangle (respectively the right-hand
side of (\BO)) for {\it arbitrary} $m$ is an infinite problem because we have
to compute M\"obius functions for $NC^m(\Phi)$ (respectively for
$(NC^m(\Phi))^*$) for $m=1,2,\dots$, this is not really true. We should
recall from (\CH) that an element of $NC^m(\Phi)$ has the form
$$\multline (w_0;w_1,\dots,w_m),\text{ with }w_0w_1\cdots w_m=c\text{ and }\\
\ell_T(w_0)+\ell_T(w_1)+\dots+\ell_T(w_m)=\ell_T(c)=n,
\endmultline\tag\AE$$
where $n$ is the rank of the root system $\Phi$. Now, $n$ can be at
most $8$ for an exceptional root system (with equality only for
$\Phi=E_8$). This implies that only at most $8$ of the $w_i$'s can be
different from the identity element $\ep$ in $W=W(\Phi)$. Hence, a
typical interval in $(NC^m(\Phi))^*$ looks like $[ u, w]$, where
$$ u=(u_0;u_1,\dots,u_m),\quad  w=(w_0;w_1,\dots,w_m),$$
$u_i=w_i=\ep$ for all but at most $8$ indices $i\ge1$,
and $u_i\le w_i$ for these remaining indices.
Let these latter indices be $i_1,i_2,\dots,i_d$, with $d\le 8$.
Then, such an interval $[ u, w]$ is
isomorphic to the ``compressed" interval $[ u', w']$, where
$$ u'=(u_0;u_{i_1},\dots,u_{i_d}),\quad 
 w'=(w_0;w_{i_1},\dots,w_{i_d}).
$$
Note that ``compressed" means that all of
$w_{i_1},w_{i_2},\dots,w_{i_d}$ are different from $\ep$.

So, what we have to do is to determine all different
{\it compressed\/} intervals $[ u', w']$. The contribution
of a compressed interval $[ u', w']$ to the right-hand side
of (\BO) is then
$$\binom md\cdot\mu^*( u', w')\,x^{\rk^* w'}y^{\rk^* u'},\tag\aE$$
because there are $\binom md$ different ways to choose
$\{i_1,i_2,\dots,i_d\}$ out of $\{1,2,\dots,m\}$. To obtain the
$M$-triangle we ``just" have to collect all these 
contributions and sum them over all possible compressed intervals.
Note that this is now a {\it finite} problem because the number of
{\it compressed\/} intervals is finite.

Rather than running through all compressed intervals, a more efficient
way to implement this is as follows. 
We rewrite the right-hand side of (\BO) as
$$ \sum _{u,w\in(NC^m(\Phi))^*}
^{}\mu^*(u,w)\,(-x)^{\rk^*w}(-y)^{\rk^*u}=
 \sum _{w\in(NC^m(\Phi))^*}
^{}(-x)^{\rk^*w}\cdot\chi^*_{[\hat 0,w]}(-y),
\tag\Ae$$
where $\hat0=(c;\ep,\dots,\ep)$ is the minimum in $(NC^m(\Phi))^*$, and
where
$$\chi^*_{[\hat 0,w]}(y)=\sum _{u\in(NC^m(\Phi))^*}
\mu^*(u,w)\,y^{\rk^*u}$$
is, essentially, the characteristic polynomial of the interval
$[\hat 0,w]$. (To be precise, it is the characteristic polynomial of
the interval $[w,\hat 0]$ in $NC^m(\Phi)$, see \cite{\StanAP, Sec.~3.10}.)
If $w=(w_0;w_1,\dots,w_m)$ with $w_{i_1},w_{i_2},\dots,w_{i_d}$
those among $w_1,w_2,\dots,w_m$ which are different from the identity
element $\ep$, then
$$[\hat 0,w]\cong[\ep,w_{i_1}]\times[\ep,w_{i_2}]\times
\dots\times[\ep,w_{i_d}],$$
where each interval $[\ep,w_{i_j}]$ is an interval in $NC(\Phi)$.
Since the characteristic polynomial is multiplicative, this implies
$$\chi^*_{[\hat 0,w]}(y)=\chi^*_{[\ep,w_{i_1}]}(y)\chi^*_{[\ep,w_{i_2}]}(y)\cdots
\chi^*_{[\ep,w_{i_d}]}(y),$$
where $\chi^*_{[\ep,w_{i_j}]}(y)=
\sum _{v\in NC(\Phi)} ^{}\mu(v,w_{i_j})y^{\rk v}$, 
with $\mu$ the M\"obius
function and $\rk$ the rank function in $NC(\Phi)$.

According to a result by Bessis \cite{\BesDAA, Lemma~1.4.3, Cor.~1.6.2}, 
each element $w_{i_j}$ is some parabolic Coxeter element (that is, a Coxeter
element in some parabolic subgroup), and the interval
$[\ep,w_{i_j}]$ is isomorphic to some $NC(\Psi)$, where $\Psi$ is the
root system of this parabolic subgroup. 

If we put all this together, then (\Ae) becomes
$$\multline
\sum _{d=0} ^{n}
\sum _{(T_1,\dots,T_d)} ^{}(-x)^{\rk T_1+\dots+\rk T_d}\cdot
N_\Phi(T_1,T_2,\dots,T_d)\\
\cdot\chi^*_{NC(T_1)}(-y)
\chi^*_{NC(T_2)}(-y)\cdots \chi^*_{NC(T_d)}(-y)\binom md,
\endmultline\tag\Aee$$
where the inner sum is over all possible $d$-tuples
$(T_1,T_2,\dots,T_d)$ of types (not necessarily irreducible types), 
and where $N_\Phi(T_1,T_2,\dots,T_d)$
is the number of ``minimal" products $c_1c_2\cdots c_d$
less than or equal to the Coxeter
element $c$ in absolute order, 
``minimal" meaning that all the $c_i$'s are different from
$\ep$ and that
$\ell_T(c_1)+\ell_T(c_2)+\dots+\ell_T(c_d)=\ell_T(c_1c_2\cdots c_d)$, 
such that the
type of $c_i$ as a parabolic Coxeter element is $T_i$, $i=1,2,\dots,d$.
The notation $NC(T)$ in (\Aee) 
means $NC(\Psi)$, where $\Psi$ is a root system
of type $T$, and $\rk T$ denotes the rank of $\Psi$.
We point out that the appearance of the binomial coefficient $\binom
md$ is explained by (\aE).

So, what we have to do to apply formula (\Aee) to compute the
right-hand side of (\BO) is, first, to determine all the 
``decomposition numbers"
$N_\Phi(T_1,T_2,\dots,T_d)$. 
Since we shall refer to it later, we point out that
these decomposition numbers have many relations
between themselves.
For example, the number $N_\Phi(T_1,T_2,\dots,T_d)$ is independent
of the order of the types $T_1,T_2,\dots,T_d$, that is, we have
$$
N_\Phi(T_{\si(1)},T_{\si(2)},\dots,T_{\si(d)})=N_\Phi(T_1,T_2,\dots,T_d)
\tag\Aa$$
for any permutation $\si$ of $\{1,2,\dots,d\}$. This follows from the
(proof of the) Shifting Lemma (see \cite{\ArmDAA, Lemma 1.3.1}). Furthermore,
by the definition of these numbers, those of ``lower rank" can be
computed from those of ``full rank." To be precise, we have
$$
N_\Phi(T_1,T_2,\dots,T_d)=
\sum _{T} ^{}N_\Phi(T_1,T_2,\dots,T_d,T),
\tag\Ab$$
where the sum is over all types $T$ of rank $n-\rk T_1-\rk T_2-
\dots -\rk T_d$ (with $n$ still denoting the rank of the 
fixed root system $\Phi$).

Second, one needs a list of the characteristic
polynomials $\chi^*_{NC(\Psi)}(y)$ for all 
{\it irreducible} root systems $\Psi$. 
(By the multiplicativity of the characteristic polynomial,
this then gives also formulae for the characteristic polynomials
of all the reducible types.)
In fact, the numbers $N_\Psi(T_1,T_2,\dots,T_d)$ carry all the
information which is necessary to do this recursively. 
Namely, by the definition of $NC(\Psi)$ and of
the decomposition numbers
$N_\Psi(T_1,T_2,\dots,T_d)$, we have
$$\chi^*_{NC(\Psi)}(y)=
\sum _{T_1,T_2} ^{}N_{\Psi}(T_1,T_2)\,
\mu_{NC(T_2)}\!\left(\hat 0_{NC(T_2)},\hat 1_{NC(T_2)}\right)
y^{\rk T_1},\tag\Ac$$
where $\mu_{NC(T_2)}(.,.)$ denotes the M\"obius function in
$NC(T_2)$,
and where $\hat 0_{NC(T_2)}$ and $\hat 1_{NC(T_2)}$ are, 
respectively, the minimal and the maximal element in $NC(T_2)$.
Indeed, inductively, the M\"obius functions 
$\mu_{NC(T_2)}\!\left(\hat 0_{NC(T_2)},\hat 1_{NC(T_2)}\right)$
are already known for all $T_2$ of lower rank than the rank
of $\Psi$. Hence, the only unknown in (\Ac) is
$\mu_{NC(\Psi)}\!\left(\hat 0_{NC(\Psi)},\hat 1_{NC(\Psi)}\right)$.
However, the latter can be computed by setting $y=1$ in (\Ac) and
using the fact that $\chi^*_{NC(\Psi)}(1)=0$ for all root
systems $\Psi$ of rank at least $1$. (This fact is
equivalent to the statement that 
$\sum _{u\in NC(\Psi)} ^{}\mu_{NC(\Psi)}\!
\left(u,\hat 1_{NC(\Psi)}\right)=0$, 
which is nothing but a part of the
definition of the M\"obius function.
Alternatively, one may use the uniform formula for the zeta 
polynomial of the non-crossing partition lattices, in which one
specializes the variable to $-1$. See 
\cite{\ChaFAA, Prop.~9}; the reader may be warned that
a slightly different convention for the zeta polynomial is
used there.)

We show in Sections~13--17 how to implement this procedure for the
dihedral root system $I_2(a)$, for the hyperbolic root systems
$H_3$ and $H_4$, and for $F_4$ and $E_6$. 
We list the values of the characteristic polynomials of the irreducible
root systems that we need below.
$$\allowdisplaybreaks\align \chi^*_{A_1}(y)&=y-1,\\
\chi^*_{A_2}(y)&=y^2-3y+2,\\
\chi^*_{I_2(a)}(y)&=y^2-ay+a-1,\\
\chi^*_{A_3}(y)&=y^3-6y^2+10y-5,\\
\chi^*_{A_4}(y)&=y^4-10y^3+30y^2-35y+14,\\
\chi^*_{A_5}(y)&=y^5-15y^4+70y^3-140y^2+126y-42,\\
\chi^*_{B_3}(y)&=y^3-9y^2+18y-10,\\
\chi^*_{D_4}(y)&=y^4-12y^3+39y^2-48y+20,\\
\chi^*_{D_5}(y)&=y^5-20y^4+106y^3-230y^2+220y-77,\\
\chi^*_{H_3}(y)&=y^3-15y^2+35y-21,\\
\chi^*_{F_4}(y)&=y^4-24 y^3+101 y^2-144 y+66,\\
\chi^*_{H_4}(y)&=y^4-60 y^3+307 y^2-480 y+232,\\
\chi^*_{E_6}(y)&=y^6-36y^5+300y^4-1035y^3+1720y^2-1368y+418.
\tag\Aeee\endalign$$


\subhead 13. Proof of the $F=M$ Conjecture for $I_2(a)$\endsubhead
By (\CA), we have
$$\multline 
(1-xy)^2F^m_{I_2(a)}\(\frac {x(1+y)} {1-xy},\frac {xy} {1-xy}\)= 
\frac{m (am-a+2)}{2}  x^2 y^2 +a m^2  x^2y\\
+\frac{m (am +a-2)}{2}  x^2+a m xy+a m 
    x+1
\endmultline\tag\BY$$
for the left-hand side of (\BO).

We now compute the right-hand side of (\BO) following the proposed
procedure in the previous section. 
We have $N_{I_2(a)}(I_2(a))=1$, $N_{I_2(a)}(A_1,A_1)=a$, 
$N_{I_2(a)}(A_1)=a$, $N_{I_2(a)}(\emptyset)=1$,
all other numbers $N_{I_2(a)}(T_1,\dots,T_d)$ being zero. Thus,
according to (\Aee) and (\Aeee), the right-hand side of (\BO) is equal
to
$$(-x)^2(y^2+ay+a-1)m+(-x)^2a(-y-1)^2\binom m2+(-x)a(-y-1)m+1,$$
which agrees with (\BY).\quad \quad \qed

\subhead 14. Proof of the $F=M$ Conjecture for $H_3$\endsubhead
By (\CB), we have
$$\multline 
(1-xy)^3F^m_{H_3}\(\frac {x(1+y)} {1-xy},\frac {xy} {1-xy}\)= 
 \frac{ m (5 m-4) (5 m-2)}{3} x^3 y^3 +5 m^2 (5 m-2) x^3 y^2 \\
+\frac{m
    (5 m+2) (5 m+4)}{3}  x^3+5 m^2 (5 m+2) x^3y+5 m (5 m-2) x^2 y^2 \\
+50 m^2 x^2y+5 m (5 m+2)
    x^2+15 m xy+15 m x+1
\endmultline\tag\Af$$
for the left-hand side of (\BO).

We now compute the right-hand side of (\BO) following the proposed
procedure in Section~12. The conclusions which we report here are
based on {\sl Maple} computations which we performed using
Stembridge's {\tt coxeter} package \cite{\StemAZ}.

We have $N_{H_3}(H_3)=1$, $N_{H_3}(A_1^2,A_1)=5$,
$N_{H_3}(A_2,A_1)=5$, $N_{H_3}(I_2(5),A_1)=5$,
$N_{H_3}(A_1,A_1,A_1)=50$, plus the assignments implied by
(\Aa) and (\Ab),
all other numbers $N_{H_3}(T_1,\dots,T_d)$ being zero. Thus,
according to (\Aee) and (\Aeee), the right-hand side of (\BO) is equal
to
$$\multline 
(-x)^3(-y^3-15y^2-35y-21)m+
2\cdot(-x)^35(-y-1)^3\binom m2\\+
2\cdot(-x)^35(y^2+3y+2)(-y-1)\binom m2+
2\cdot(-x)^35(y^2+5y+4)(-y-1)\binom m2\\+
(-x)^350(-y-1)^3\binom m3+
(-x)^25(-y-1)^2m+
(-x)^25(y^2+3y+2)m\\+
(-x)^25(y^2+5y+4)m+
(-x)^250(-y-1)^2\binom m2+
(-x)15(-y-1)m+1,
\endmultline$$
which agrees with (\Af).\quad \quad \qed

\subhead 15. Proof of the $F=M$ Conjecture for $H_4$\endsubhead
By (\CC), we have
$$\allowdisplaybreaks\multline 
(1-xy)^4F^m_{H_4}\(\frac {x(1+y)} {1-xy},\frac {xy} {1-xy}\)= 
 \frac{1}{4} m (3 m-1) (5 m-3) (15 m-14) x^4 y^4 \\
+15 m^2 (3 m-1) (5 m-3) x^4 y^3
    +\frac{1}{2} m^2 \left(675 m^2-61\right) x^4 y^2 \\
+15 m^2 (3 m+1) (5 m+3) x^4y+\frac{1}{4} m (3 m+1)
    (5 m+3) (15 m+14) x^4\\
+15 m (3 m-1) (5 m-3)
    x^3 y^3 
+15 m^2 (45 m-14) x^3 y^2 \\
+15 m^2 (45 m+14)
    x^3y+15 m (3 m+1) (5 m+3) x^3\\
+\frac{1}{2} m (465 m-149) x^2 y^2 +465
    m^2 x^2y+\frac{1}{2} m (465 m+149) x^2+60 m xy+60 m x+1
\endmultline\tag\AF$$
for the left-hand side of (\BO).

We now compute the right-hand side of (\BO) following the proposed
procedure in Section~12. The conclusions which we report here are
based on {\sl Maple} computations which we performed using
Stembridge's {\tt coxeter} package \cite{\StemAZ}.

We have $N_{H_4}(H_4)=1$, $N_{H_4}(A_1*A_2,A_1)=15$,
$N_{H_4}(A_3,A_1)=15$, 
$N_{H_4}(H_3,A_1)=15$,
$N_{H_4}(A_1*I_2(5),A_1)=15$,
$N_{H_4}(A_1^2,A_1^2)=30$,
$N_{H_4}(A_1^2,A_2)=30$,
$N_{H_4}(A_1^2,I_2(5))=15$,
$N_{H_4}(A_2,A_2)=5$,
$N_{H_4}(A_2,I_2(5))=15$,
$N_{H_4}(I_2(5),I_2(5))=3$,
$N_{H_4}(A_1^2,A_1,A_1)=225$,
$N_{H_4}(A_2,A_1,A_1)=150$,
$N_{H_4}(I_2(5),A_1,A_1)=90$,
$N_{H_4}(A_1,A_1,A_1,A_1)=1350$, 
plus the assignments implied by
(\Aa) and (\Ab),
all other numbers $N_{H_4}(T_1,\dots,T_d)$ being zero. 
If one substitutes accordingly in (\Aee), using the information from
(\Aeee), then one obtains an expression which agrees with (\AF)
after simplification.\quad
\quad \qed

\subhead 16. Proof of the $F=M$ Conjecture for $F_4$\endsubhead
By (\CD), we have
$$\multline 
(1-xy)^4F^m_{F_4}\(\frac {x(1+y)} {1-xy},\frac {xy} {1-xy}\)= 
 \frac{1}{2} m (2 m-1) (3 m-1) (6 m-5)  x^4y^4\\
+12 m^2 (2
m-1) (3 m-1) 
    x^4y^3 +m^2 \left(108 m^2-7\right) x^4y^2\\
+12 m^2 (2 m+1) (3 m+1) x^4y+\frac{1}{2} m (2 m+1) (3 m+1) (6
    m+5) x^4\\
+12 m (2 m-1) (3 m-1)  x^3y^3+12 m^2
    (18 m-5)  x^3y^2\\
+12 m^2 (18 m+5) x^3y+12 m (2 m+1) (3 m+1) x^3
+m (78
    m-23)  x^2y^2\\+156 m^2 x^2y+m (78 m+23) x^2
+24 m xy+24 m x+1 
\endmultline\tag\AG$$
for the left-hand side of (\BO).

We now compute the right-hand side of (\BO) following the proposed
procedure in Section~12. The conclusions which we report here are
based on {\sl Maple} computations which we performed using
Stembridge's {\tt coxeter} package \cite{\StemAZ}.

We have $N_{F_4}(F_4)=1$, $N_{F_4}(A_1*A_2,A_1)=12$,
$N_{F_4}(B_3,A_1)=12$, 
$N_{F_4}(A_1^2,A_1^2)=12$,
$N_{F_4}(A_1^2,B_2)=12$,
$N_{F_4}(A_2,A_2)=16$,
$N_{F_4}(B_2,B_2)=3$,
$N_{F_4}(A_1^2,A_1,A_1)=72$,
$N_{F_4}(A_2,A_1,A_1)=48$,
$N_{F_4}(B_2,A_1,A_1)=36$,
$N_{F_4}(A_1,A_1,A_1,A_1)=432$, 
plus the assignments implied by
(\Aa) and (\Ab),
all other numbers $N_{F_4}(T_1,\dots,T_d)$ being zero. 
If one substitutes accordingly in (\Aee), using the information from
(\Aeee), then one obtains an expression which agrees with (\AG)
after simplification.\quad
\quad \qed

\subhead 17. Proof of the $F=M$ Conjecture for $E_6$\endsubhead
By (\CE), we have
{\eightpoint
$$\multline 
(1-xy)^6F^m_{E_6}\(\frac {x(1+y)} {1-xy},\frac {xy} {1-xy}\)= 
 \frac{1}{30} m (2 m-1) (3 m-1) (4 m-1) (6 m-5) (12 m-7)  x^6y^6\\
+\frac{6}{5}
    m^2 (2 m-1) (3 m-1) (4 m-1) (12 m-7)  x^6y^5
+2 m^2 (3 m-1) (4 m-1)
    \left(36 m^2-9 m-2\right)  x^6y^4\\
+3 m^2 (4 m-1) (4 m+1) \left(24
    m^2-1\right)  x^6y^3+2 m^2 (3 m+1) (4 m+1) \left(36 m^2+9 m-2\right) 
    x^6y^2\\
+\frac{6}{5} m^2 (2 m+1) (3 m+1) (4 m+1) (12 m+7)  x^6y
+\frac{1}{30} m (2 m+1) (3 m+1) (4 m+1) (6 m+5) (12 m+7)
    x^6\\
+\frac{6}{5} m
    (2 m-1) (3 m-1) (4 m-1) (12 m-7)  x^5y^5
+12 m^2 (3 m-1) (4 m-1) (12 m-5)
     x^5y^4\\
+6 m^2 (4 m-1) \left(144 m^2-12 m-5\right)  x^5y^3+6 m^2 (4 m+1)
    \left(144 m^2+12 m-5\right)  x^5y^2\\
+12 m^2 (3 m+1) (4 m+1) (12 m+5)  x^5y
+\frac{6}{5} m (2 m+1) (3 m+1) (4 m+1)
    (12 m+7) x^5\\
+2 m (3 m-1) (4 m-1) (30
    m-13)  x^4y^4+6 m^2 (4 m-1) (120 m-31)  x^4y^3+16 m^2 \left(270
    m^2-7\right)  x^4y^2\\
+6 m^2 (4 m+1) (120
    m+31)  x^4y
+2 m (3 m+1) (4 m+1) (30 m+13) x^4
+9 m (4 m-1) (18 m-5)  x^3y^3\\
+18 m^2 (108 m-19)  x^3y^2
+18 m^2 (108 m+19)  x^3y+9 m (4
    m+1) (18 m+5) x^3+12 m (21 m-4)  x^2y^2\\
+504 m^2  x^2y+12 m (21
    m+4) x^2+36 m xy+36 m x+1\\
\endmultline\tag{\tenpoint\AH}$$}%
for the left-hand side of (\BO).

We now compute the right-hand side of (\BO) following the proposed
procedure in Section~12. The conclusions which we report here are
based on {\sl Maple} computations which we performed using
Stembridge's {\tt coxeter} package \cite{\StemAZ}.

We have $N_{E_6}(E_6)=1$, 
 $N_{E_6}(A_1*A_2^2, A_1) = 6$,
 $N_{E_6}(A_1*A_4, A_1) = 12$, $N_{E_6}(A_5, A_1) = 6$,  
 $N_{E_6}(D_5, A_1) = 12$, 
 $N_{E_6}(A_1^2*A_2, A_2) = 36$, $N_{E_6}(A_2^2, A_2) = 8$, 
 $N_{E_6}(A_1*A_3, A_2) = 24$, $N_{E_6}(A_4, A_2) = 24$, $N_{E_6}(D_4, A_2) = 4$, 
 $N_{E_6}(A_1^2*A_2, A_1^2) = 18$, 
 $N_{E_6}(A_1*A_3, A_1^2) = 36$, $N_{E_6}(A_4, A_1^2) = 36$, $N_{E_6}(D_4, A_1^2) = 18$, 
 $N_{E_6}(A_1^3, A_1^3) = 12$, $N_{E_6}(A_1*A_2, A_1^3) = 24$,\linebreak $N_{E_6}(A_1*A_2, A_1*A_2) = 48$, 
 $N_{E_6}(A_3, A_1^3) = 36$, $N_{E_6}(A_3, A_1*A_2) = 72$, $N_{E_6}(A_3, A_3) = 27$, 
 $N_{E_6}(A_1^2*A_2, A_1, A_1) = 144$, $N_{E_6}(A_2^2, A_1, A_1) = 24$, 
 $N_{E_6}(A_1*A_3, A_1, A_1) = 144$, $N_{E_6}(A_4, A_1, A_1) = 144$, $N_{E_6}(D_4, A_1, A_1) = 48$, 
 $N_{E_6}(A_1^3, A_1^2, A_1) = 180$, $N_{E_6}(A_1^3, A_2, A_1) = 168$, 
 $N_{E_6}(A_1*A_2, A_1^2, A_1) = 360$,
 $N_{E_6}(A_1*A_2, A_2, A_1) = 336$, 
 $N_{E_6}(A_3, A_1^2, A_1) = 378$, $N_{E_6}(A_3, A_2, A_1) = 180$, 
 $N_{E_6}(A_1^2, A_1^2, A_1^2) = 432$, $N_{E_6}(A_2, A_1^2, A_1^2) = 504$, 
 $N_{E_6}(A_2, A_2,\mathbreak A_1^2) = 288$, $N_{E_6}(A_2, A_2, A_2) = 160$,
 $N_{E_6}(A_1^2, A_1^2, A_1, A_1) = 2376$, 
 $N_{E_6}(A_2, A_1^2, A_1, A_1) = 1872$, $N_{E_6}(A_2, A_2, A_1, A_1) = 1056$,
 $N_{E_6}(A_1^3, A_1, A_1, A_1) = 864$, $N_{E_6}(A_1*A_2, A_1, A_1, A_1) = 1728$, 
 $N_{E_6}(A_3, A_1, A_1, A_1) = 1296$, 
 $N_{E_6}(A_1^2, A_1, A_1, A_1, A_1) = 10368$, $N_{E_6}(A_2, A_1, A_1, A_1,\mathbreak A_1) = 6912$,
 $N_{E_6}(A_1, A_1, A_1, A_1, A_1, A_1) = 41472$,
plus the assignments implied by
(\Aa) and (\Ab),
all other numbers $N_{E_6}(T_1,\dots,T_d)$ being zero. 
If one substitutes accordingly in (\Aee), using the information from
(\Aeee), then one obtains an expression which agrees with (\AH)
after simplification.\quad
\quad \qed

\subhead 18. The dual $F$-triangle\endsubhead
Armstrong \cite{\ArmDAA, Sec.~4} defines the {\it dual $F$-triangle},
denoted here by $\tilde F^m_\Phi(x,y)$, as
$$\tilde F^m_\Phi(x,y)=(-1)^nF^m_\Phi(-1-x,-1-y),$$
where $n$ is the rank of the root system $\Phi$. He conjectures that
the dual $F$-triangle can be expressed in form of a weighted bivariate
generating function for the faces of $\De^m(\Phi)$ involving the
{\it Fuss--Narayana numbers} $\Nar^m(\Phi,i)$, the latter enumerating
all elements of rank $i$ in the $m$-divisible non-crossing partition
poset $NC^m(\Phi)$. For explicit formulae for the Fuss--Narayana numbers see
\cite{\ArmDAA, Theorem~2.3.6}. These numbers occur also as $h$-numbers
in \cite{\FoReAA, Theorem~9.2}. (One has to reverse the ordering of
the numbers to convert one sequence of numbers into the other.) In view
of our proof below, Armstrong's conjecture becomes the following theorem.

\proclaim{Theorem DF} For any finite root system $\Phi$, we have
$$\tilde F^m_\Phi(x,y)=
\sum _{k,l\ge0} ^{}\frac {\Nar^m(\Phi,k+l)}
{\Nar^1(\Phi,k+l)}f_{k,l}\, x^ky^l.
\tag\BZ$$
\endproclaim

\demo{Proof}
Clearly, for the exceptional root systems one can verify (\BZ)
routinely by using the explicit formulae for the refined face numbers,
as given through the formulae for the $F$-triangle in Section~7, and
the formulae for the Fuss--Narayana numbers in \cite{\ArmDAA,
\FoReAA}.

To verify (\BZ) for the root systems $A_n$, $B_n$ and $D_n$, some work
has to be done. However, the verifications in these types are very
similar to each other so that we give below only the proof in
type $A_n$, leaving the proofs for $B_n$ and $D_n$ to the reader.

By Theorem~FA, in type $A_n$ the left-hand side of (\BZ) is equal to
$$\allowdisplaybreaks\align 
(-1)&^n\sum _{k,l\ge0} ^{}\frac{l+1 }{k+l+1}{\binom {n} { k+l}} 
{\binom {m (n+1)+k-1} k}(-1-x)^k(-1-y)^l\\
&=\sum _{k,l,r,s\ge0} ^{}\frac{l+1 }{k+l+1}{\binom {n} { k+l}} 
{\binom {m (n+1)+k-1} k}\binom kr\binom ls(-1)^{n+k+l}x^ry^s\\
&=\sum _{k,r,s\ge0} ^{}\frac{s+1 }{n+1}{\binom {n-s-1} {n-k-s}} 
{\binom {m (n+1)+k-1} k}\binom kr(-1)^{n+k+s}x^ry^s\\
&=\sum _{k,r,s\ge0} ^{}\frac{s+1 }{n+1}{\binom {n-s-1} {n-k-s}} 
{\binom {m (n+1)+k-1} {k-r}}\binom {m(n+1)+r-1}r(-1)^{n+k+s}x^ry^s\\
&=\sum _{r,s\ge0} ^{}\frac{s+1 }{n+1}
{\binom {m (n+1)} {n-s-r}}\binom {m(n+1)+r-1}rx^ry^s.
\endalign$$
As earlier, for the evaluation of the sums over $l$ and $k$ we used
special instances of the Chu--Vandermonde summation. 

On the other hand, by Theorem~FA and by \cite{\ArmDAA, Theorem~2.3.6},
the right-hand side of (\BZ) is equal to
$$
\sum _{k,l\ge0} ^{}\frac {\frac {1} {n+1}\binom {n+1}{k+l}\binom
{m(n+1)}{n-k-l}} {\frac {1} {n+1}\binom {n+1}{k+l}\binom
{n+1}{n-k-l}}
\frac{l+1 }{k+l+1}{\binom {n} { k+l}} 
{\binom {m (n+1)+k-1} k}x^ky^l,$$
which is exactly the same expression.\quad \quad \qed
\enddemo

\bigskip
\remark{Acknowledgement} 
I am indebted to Sergey Fomin for suggesting to me to embark 
on this project and,
thus, to the Institut Mittag--Leffler, Anders Bj\"orner and Richard
Stanley for hosting both of
us at the same time at the Institut during
the ``Algebraic Combinatorics" programme in Spring 2005.
Moreover, I wish to thank
Drew Armstrong, Christos Athanasiadis and Eleni Tzanaki
for the many extremely helpful exchanges that I had with them
on the subject matter, and for making drafts for \cite{\ArmDAA},
\cite{\AthaAI} and \cite{\TzanAA} available to me at an early stage.
Last, but not least, I am indebted to the referee for
an extremely careful reading of the manuscript.
\endremark

\enddocument